\documentclass[11pt,a4paper]{article}

\usepackage{epsf}
\usepackage{psfrag}
\usepackage{amsmath,amssymb,cite}
\usepackage{lscape}
\usepackage{latexsym}
\usepackage[usenames]{color}
\usepackage[mathcal]{eucal}
\usepackage{stmaryrd}
\usepackage{textcomp}
\usepackage{setspace}  
\usepackage{graphicx}   
\usepackage{pgfplots}    
\usepackage{epsfig}       
\usepackage{epstopdf}   
\usepackage{appendix}
\usepackage{afterpage}
\usepackage{authblk}  
\usepackage[runin]{abstract}	

\newtheorem{theorem}{Theorem}[section]

\newtheorem{proposition}{Proposition}[section]

\newtheorem{remark}{Remark}[section]

\newtheorem{definition}{Definition}[section]

\usepackage{hyperref}
\hypersetup{
    colorlinks=true,
    linkcolor=blue,
    filecolor=magenta,      
    urlcolor=cyan,
    citecolor=blue,
}
\newcommand{\ds}{\displaystyle}
\newcommand{\dZ}{{\cal Z \kern -0.7em Z}}
\newcommand{\dC}{{\rm\hbox{C \kern-0.8em\raise0.2ex\hbox{\vrule height5.4pt width0.7pt}}}}
\newcommand{\dQ}{{\rm\hbox{Q \kern-0.85em\raise0.25ex\hbox{\vrule height5.4pt width0.7pt}}}}

\newcommand{\R}{\mathbb{R}}

\newcommand{\proofbox}{\hspace{\fill}{$\Box$}}

\newenvironment{proof}{Proof.}{\proofbox}

\newcommand\old[1]{}
\newcommand{\beqa}{\begin{eqnarray*}}
\newcommand{\eeqa}{\end{eqnarray*}}
\newcommand*\samethanks[1][\value{footnote}]{\footnotemark[#1]}

\begin{document}

\title{\bf Algorithms for Generating Pareto fronts of Multi-objective Integer and Mixed-Integer Programming Problems}

\author{R. S. Burachik\thanks{Mathematics, UniSA STEM, University of South Australia, Australia;\\ \{regina.burachik\} or \{yalcin.kaya\}\! @unisa.edu.au\,, or\ rizmm001@mymail.unisa.edu.au (*corresponding author)\,.}
\qquad C. Y. Kaya\samethanks
\qquad M. Mustafa Rizvi\samethanks\ \thanks{Department of Mathematics, University of Chittagong, Bangladesh.}}

\maketitle

\vspace*{-5mm}
\begin{abstract} {\sf Multi-objective integer or mixed-integer programming problems typically have disconnected feasible domains, making the task of constructing an approximation of the Pareto front challenging. The present paper shows that certain algorithms which were originally devised for continuous problems can be successfully adapted to approximate the Pareto front for integer, and mixed-integer, multi-objective optimization problems.  Relationships amongst various scalarization techniques are established to motivate the choice of a particular scalarization in these algorithms. The proposed algorithms are tested by means of two-, three- and four-objective integer and mixed-integer problems, and comparisons are made. In particular, a new four-objective algorithm is used to solve a rocket injector design problem with a discrete variable, which is a challenging mixed-integer programming problem.}
\end{abstract}

\begin{verse}
\end{verse}
\begin{verse} 
\end{verse}

\pagestyle{myheadings}
\markboth{}{\sf\scriptsize Algorithms for Multi-objective Mixed-Integer Programming \ \ by R. S. Burachik, C. Y. Kaya and M. M. Rizvi}

\section{Introduction}\label{sec1}
Multi-objective optimization is concerned with simultaneous minimization of several conflicting objective functions. The aim is to obtain/approximate, and view, the set of trade-off, or compromise, solutions.  This set of trade-off solutions is referred to, in the current article, as the ({\em weak}) {\em Pareto front}, or the ({\em weak}) {\em efficient set}, of the problem---a more precise definition is to be provided in Section~\ref{sec2}.  Once the set of Pareto points is at hand, the decision maker conveniently selects a suitable solution from this set, usually based on an additional criterion.  These kinds of problems arise in many multidisciplinary applications, for example, in health \cite{SunDepEva2014,CabEhrMasPhi2014,ZhaZha1999}, engineering \cite{MarAro2004,StaMatSta2012}, mining \cite{YuZheGaoYan2017,KauMukBas2015}, and finance \cite{Romanko2010,Smimou2014}.

If all variables of a multi-objective optimization problem belong to a continuous set, then the problem is called a {\em multi-objective continuous optimization} problem.  If, on the other hand, the variables belong to a subset of integers, then the problem is a {\em multi-objective integer programming} one.  When the model involves both integer and continuous variables the problem is referred to as a {\em multi-objective mixed-integer programming} problem.  The latter case arises in many real-life applications, for example, the knapsack \cite{KirSay2014,CeyKokLok2019} and shortest-path problems \cite{CurRevCoh1985,UluTeg1994}.  Other types of mixed-integer multi-objective problems can be found in \cite{RaiAlvCar2014}.

It is desirable that algorithms for solving multi-objective integer and mixed-integer programming problems will have the following major practical attributes: (i) construct a faithful approximation of the Pareto front and (ii) generate this approximation in a reasonable amount of computational time. Attribute (ii) is particularly valuable, since, for example, solving each scalarized problem for a mixed-integer programming problem can be very costly.

Indeed it is in general not realistic to compute too many points in the Pareto front, especially if the Pareto front contains connected segments (e.g. segments of curves in the case of two objectives and segments of surfaces in the case of three objectives).  Therefore we aim to construct an {\em approximation} of the Pareto front in the following sense.  We consider some discretization grid, or a finite partition, of the space of all weights for a parametric scalarization of the multi-objective problem.  

In Sections~\ref{sec3.1}--\ref{sec3.3}, relationships amongst various scalarization techniques are presented. Namely, it is shown that every point in the Pareto front that can be generated by the weighted-sum~\cite{Gass1955, Zadeh1963}, $k$th-objective $\epsilon$-constraint~\cite{ChaHai1983} or feasible-value-constraint scalarization~\cite{BurKayRiz2017} technique, can also be generated by the weighted-constraint scalarization technique introduced in~\cite{BurKayRiz2017} and implemented in various algorithms in \cite{BurKayRiz2017, BurKayRiz2014}.  On the other hand, a solution of the weighted-constraint scalarization may not be a solution of the other three scalarization techniques. Therefore, in the numerical experiments, the weighted-sum, $k$th objective $\epsilon$-constraint and feasible-value-constraint scalarization techniques are not incorporated into the algorithms that are used to approximate the Pareto front.

As a result, the scalarization techniques we use assure that the solution of the scalarized problem with each weight in the grid yields a Pareto solution, and the collection of all these solutions constitutes a finite and faithful subset of the Pareto front.  This finite subset is an approximation  generated by the weak efficient points; in other words, all the points in this finite set are non-dominated by each other. 

The phrase, {\em approximation of the Pareto front}, is used in the sense that, as the discretization is made finer and finer, a larger subset of the true Pareto front is generated.  The same notion of approximation was used in the authors' previous studies in \cite{BurKayRiz2014, BurKayRiz2017}.

There is an abundance of practical algorithms for constructing an approximation of the Pareto front of continuous-variable problems, especially those with two or three objective functions~\cite{BurKayRiz2014, BurKayRiz2017, MueGraeSch2009, DasDen1998, DutKay2011, Eichfelder2008, KayMau2014, ParZilZil2017, Rizvi2013, Yang1997, Ruzika2005}.  A successful algorithm in the case of four objectives has recently been proposed and numerically studied by the authors of the current paper in~\cite{BurKayRiz2017}. 

As for integer and mixed-integer programming problems, the question of generating a faithful approximation of the Pareto front is still a developing area for problems with three or more objectives.  
Multi-objective integer problems (MOIP), where all variables are integer, have been studied and effective algorithms proposed in earlier and recent literature, see e.g.\ \cite{Mela2007,Dachert2014,Boland2016,Santis2019}. The present paper focuses more on the mixed-integer case, acknowledging the fact that mixed-integer programming problems are in general far more challenging than the integer ones.  For problems which are posed in {\em convex form}, the reader is referred to the branch-and-bound approach presented in the recent article~\cite{Santis2019}. Those approaches that do work with general three objectives, do not attempt to construct the whole front---see \cite{AntAlvCli2016, BelSoyWie2013, Mavrotas2013, OzlBurMac2014, OzlAzi2009, PetOzl2017, PrzGanEhr2010a, PrzGanEhr2010b, NieEic2019}.  Hence there is a need for new techniques that approximate Pareto fronts for the mixed-integer case. 	Given the fact that for these type of problems there is no certificate of optimality, to solve them numerically, one needs a suitable solver for the subproblems that guarantees to yield a globally optimal solution. To the best of the authors' knowledge, such a solver is not (yet) available for mixed-integer programming problems.

Presumably, one of the main reasons for which there are almost no studied examples of three-\ and four-objective mixed-integer optimization problems in the literature is that the feasible set in the presence of integer variables is often disconnected.  In this situation, specialized algorithms are needed to construct the Pareto front.  Specialized scalarization techniques for problems with continuous variables and disconnected domains were first studied by the authors in~\cite{Rizvi2013, BurKayRiz2014} and then in~\cite{BurKayRiz2017}---also see~\cite[Sections 2.4 and 3.2.4]{ParZilZil2017}. 

Approximation of the Pareto \textcolor{blue}{front} indeed relies on the success of the scalarization technique employed to transform the problem into a single variable one. If the multi-objective problem is continuous, and the domain is connected, then the Pascoletti--Serafini scalarization \cite{PasSer1984, Eichfelder2008}, or the constrained-Tchebycheff scalarization \cite{DutKay2011}, can be successfully used to construct the Pareto front. If the domain is disconnected, however, the Pascoletti--Serafini scalarization, as well as most of the other available scalarization techniques, are likely to fail, as illustrated in \cite{BurKayRiz2017, BurKayRiz2014}.  In such cases, implementation of certain other scalarization techniques, in particular the {\em weighted-constraint scalarization} technique and the associated algorithms introduced in \cite{BurKayRiz2017, BurKayRiz2014}, have been shown to perform well in constructing the Pareto front.  

The algorithms in \cite{BurKayRiz2017, BurKayRiz2014} have been demonstrated to be particularly successful when applied to three- and four-objective {\em continuous} optimization problems.  However, the same algorithms have never been tested before on multi-objective integer or mixed-integer optimization problems.  Therefore it is natural to extend and implement the algorithms in \cite{BurKayRiz2017, BurKayRiz2014} for mixed-integer problems.  If a problem is {\em linear} then the solver PolySCIP developed by Schenker \cite{Schenker2019} or the search based algorithms introduced in \cite{BolChaSav2017}, can possibly be used.  Hence, there is a need to devise algorithms for {\em nonlinear} and {\em nonconvex} mixed-integer optimization problems, in particular for problems with more than three objective functions. Note that an integer or mixed-integer problem is called nonconvex when the objective function or the {\em form of the constraints} is not convex.

The current paper shows that Algorithms~1~\cite{BurKayRiz2014}, 3 and 5~\cite{BurKayRiz2017} are particularly successful in solving multi-objective integer or mixed-integer programming problems, compared with Algorithms 2~\cite{BurKayRiz2014}, 4~\cite{BurKayRiz2017} and 6~\cite{MueGraeSch2009, BurKayRiz2017}.  The latter group of algorithms use the conventional Pascoletti--Serafini scalarization, which is effective in the case of simply connected Pareto fronts.  Moreover, the new Algorithm~7 is an extension of Algorithm~5 to the case of four objective functions.  Algorithm~7 itself and its application to a rocket injector design problem~\cite{BurKayRiz2017, GoeVaiHafShyQueTuc2007, VaiTucPapShy2004} in Section~\ref{rocket} is the first of its kind, as, to the knowledge of the authors, there does not exist a multi-objective mixed-integer optimization algorithm applicable to problems with three or more {\em nonlinear} objective functions.  See Table~\ref{algorithms} and Section~\ref{descriptions} for further details.

The main contributions of the present paper are summarized as follows.
\begin{itemize}
\item A theoretical comparison study of the four scalarization techniques mentioned above is carried out so as to motivate the choice of the particular techniques employed in the numerical implementation of the algorithms. 
\item Algorithm~5 given in \cite{BurKayRiz2017} is extended to the case of four-objective mixed-integer programming problems, and implemented for a challenging test problem.
\item  Algorithm~1 given in \cite{BurKayRiz2014}---also see \cite[Section 3.2.4]{ParZilZil2017} and Algorithm~8.1 given in \cite{MueGraeSch2009} are tested and compared for three-objective integer programming problems. Note that these algorithms were previously tested in \cite{BurKayRiz2014} and \cite{MueGraeSch2009} only on continuous optimization problems.
\end{itemize}

The layout of the paper is as follows. The problem description along with some preliminaries are presented in Section~\ref{sec2}.  Five scalarization approaches, namely the Pascoletti--Serafini, weighted-sum, $k$th-objective $\epsilon$-constraint, feasible-value-constraint, and weighted-constraint scalarization, are recalled in Section~\ref{sec3}.  In this section, relationships among several of these scalarization techniques are established. In Section~\ref{sec4}, algorithms for two and three objectives are described, and a new algorithm for four objectives is introduced.  Numerical experiments with these algorithms are given in Section~\ref{sec5}.  In particular, Section~\ref{rocket} presents the results for a new and challenging reformulation of the rocket injector design problem involving four objective functions. Conclusion and discussion is given in Section \ref{sec6}. Detailed descriptions of Algorithms~3 and 7 are included in the Appendix.

\section{Preliminaries}
\label{sec2}

In this section, we collect the relevant notions, definitions, and concepts that are used in our study.  These are standard notations and tools for multi-objective optimization, and we follow here the classical notation in the literature---see for example, \cite{ChaHai1983, Miettinen1999,Yu1985}.  Let $\mathbb{R}^n$ be the $n$-dimensional Euclidean space. Define $\mathbb{R_+}$ to be the set of all nonnegative real numbers and $\mathbb{R_{++}}$ the set of positive numbers. Given the positive integers $\ell,m,n_1$ and $n_2$, we consider the following multi-objective optimization Problem~$(P)$.
\begin{equation}  \tag{$P$}
\left\{\begin{array}{rl} 
\ds\min & \ \ f(x) := \left[f_1(x),\ldots,f_{\ell}(x)\right] \\[2mm]
\mbox{s.t.} & \ \ x \in  X := \left\{ x \in \mathbb{Z}^{n_1} \times \mathbb{R}^{n_2} \mid g_j(x) \leq  0, \;\;j=1,\ldots,m \right\},
\end{array}
\right.
\end{equation}
where $\mathbb{Z}$ is set of all integers, and $f_i:\mathbb{Z}^{n_1} \times \mathbb{R}^{n_2} \to \mathbb{R}$, $i=1,\ldots,\ell$, and $g_j:\mathbb{Z}^{n_1} \times \mathbb{R}^{n_2} \to \mathbb{R}$, $j=1,\ldots,m$. We assume that the functions $f_i$ are bounded below on the constraint set $X$, and the lower bounds of the functions $f_i$ are known. Therefore, we can impose that
\begin{equation}\label{positive}
\ds\min_{i=1,\ldots,\ell}\ \left\{\min_{x\in  X} \;\; f_i(x)\right\} > 0\,.
\end{equation}
We will use the following set of positive weights in our analysis.
\[\hspace*{-5mm} W := \left
\{ w \in \mathbb{R}^{\ell} \mid w_i > 0\,,\ \sum_{i=1}^{\ell} w_i=1
\right \}.\]

We recall next two types of solutions of $(P)$. The more restrictive type of solution is the so-called {\em efficient point} \cite{Yu1985} or {\em Pareto point} \cite{Miettinen1999}, and a less restrictive concept is the one of a {\em weak efficient} point.   We also provide below the concepts of {\em ideal} and {\em utopia} {\em vectors} and {\em individual minima}.

\begin{definition} \label{edef}\rm(See \cite{ChaHai1983} and \cite{Miettinen1999})
\begin{itemize}
\item[(a)]  A point $\bar{x} \in X$ is said to be \textit{efficient} for Problem $(P)$ iff there is no $x \in X$, such that $f_i(x) \leq f_i(\bar{x})$, $\forall\,i\in\{1,\ldots,\ell\}$, and $f_j(x) < f_j(\bar{x})$, for some $j\in\{1,\ldots,\ell\}$. Let, $E(P)$ be the set of efficient points of Problem $(P)$.
\item[(b)]  A point $\bar{x} \in X$ is said to be \textit{weak efficient} for Problem $(P)$ iff there is no $x \in X$ such that $f_i(x) < f_i(\bar{x})$, $\forall\,i\in\{1,\ldots,\ell\}$.  Let,  $WE(P)$ be the set of weak efficient points of Problem~$(P)$.  We
define the {\em Pareto front} of Problem~$(P)$ as the image of $WE(P)$ under $f$.
\item[(c)] Suppose that $\bar{x}_{f_i}$ is a minimizer of $f_i$, $i = 1,\ldots,\ell$, over the set $X$.  That is to say, $\bar{x}_{f_i}$ solves the optimization problem
\begin{equation} \tag{P$_i$}
\min_{x\in X}\ f_i(x)\,.
\end{equation}
Then the vector $\bar{x}_f = [\bar{x}_{f_1},\ldots,\bar{x}_{f_\ell}]$ is called an \emph{ideal vector}.  We denote the vector of \emph{individual minima} (IM) by $f(\bar{x}_{f_i})=[f_1(\bar{x}_{f_1}),f_2(\bar{x}_{f_2}),\ldots, f_{\ell}(\bar{x}_{f_{\ell}})]$.
\item[(d)]\label{utopia} \rm
\textit{A utopia vector} $u := (u_1,\ldots, u_\ell) \in \mathbb{R}^\ell$ associated with Problem~$(P)$ is defined by $u_i := f_i(\bar{x}_{f_i}) - \varepsilon_i$,  where, for all $i=1,\ldots,\ell$, $\varepsilon_i > 0$.
\end{itemize}
\end{definition}

For the problem to be truly multi-objective, ideal vectors cannot be feasible. It is worth noting that $E(P)$ $\subseteq$  $WE(P)$,  but the opposite inclusions in general do not hold.

\medskip

At first glance, weak efficient solutions may not appear to be convenient from a practical point of view, since they could (in theory) be improved.  On the other hand, some Pareto fronts are extremely difficult to approximate (as in the case of the rocket injector design problem that we study in Section~\ref{rocket}), and practitioners may have to look for weak efficient solutions which are not efficient as well as the efficient ones.  In other words, in situations when finding an efficient solution is not an easy task, weak efficient solutions need to be welcome.  

Weak efficient points may also help eliminate (or filter out) some of the dominated points produced by the numerical procedure. In particular, they allow one to obtain necessary conditions for non-dominated solutions---see~\cite{GRTZ2003, Luc1989}. 

Besides from these valuable theoretical aspects of weak efficient points, we refer to the works \cite{Yu1985,Jahn2004, Stewart1997,ChaHai1983}, where the approximation of the whole Pareto Front (efficient and weak efficient points) has important applications. 

Studies of existence of weak efficient points, such as the ones in \cite{Gourion2008,Gutierrez2020} demonstrate their relevance in multi-objective optimization. A concrete example of the usefulness of approximating the entire Pareto front (efficient and weak efficient points), is illustrated later on via one of our test problems in Remark \ref{rem:5.1}.

\medskip

A comment is in order regarding the different parts of the Pareto front.  In what follows, we will be using the expressions {\em interior} and/or {\em boundary} of the Pareto front. We clarify next what we mean exactly by these latter two expressions, since they are not to be confused with the classical concepts of interior and boundary used in topology theory. 

Since the Pareto front has no interior in $\R^{\ell}$ (otherwise we could find points in the front which are strictly dominated), it is homeomorphic to a subset of $\R^{\ell-1}$. Hence, when we refer to the {\em boundary} of the Pareto front, what we actually mean is the subset of the front (in the $\R^\ell$ space), which is homeomorphic to the topological boundary of its homeomorphic image in $\R^{\ell-1}$. Similarly, by {\em interior} of the Pareto front, we mean the subset of the front (in the $\R^\ell$ space), which is homeomorphic to the topological interior of its homeomorphic image in $\R^{\ell-1}$.

For instance, if $\ell=2$, the boundary of the Pareto front (if connected) consists of the points $(\min_{x\in X}f_1(x), \max_{x\in X}f_2(x))$ and $(\min_{x\in X}f_2(x), \max_{x\in X}f_1(x))$.  If $\ell=3$ and the Pareto front in $\R^3$ is connected, then it is homeomorphic to the two-dimensional simplex. In this case, the boundary of the Pareto front is the union of regions homeomorphic to lower dimensional faces of the simplex.  To avoid the use of involved terminology, we refer in what follows, to the boundary and interior of the Pareto front via this homeomorphism.

\section{Scalarization Techniques and Their Relationships} 
\label{sec3}

In this section, we recall five scalarization techniques; namely, the
\begin{itemize}
\item[(i)] \textit{weighted-sum scalarization}~\cite{Gass1955}, 
\item[(ii)] \textit{$k$th-objective $\epsilon$-constraint scalarization}~\cite{ChaHai1983}, 
\item[\textcolor{blue}{(iii)}] \textit{Pascoletti--Serafini scalarization}~\cite{PasSer1984},  
\item[(iv)] ({\em $k$th-objective}) {\em weighted-constraint scalarization}~\cite{BurKayRiz2014},
\item[(v)] \textit{feasible-value-constraint scalarization}~\cite{BurKayRiz2017}. 
\end{itemize}
The first three of these have been popular/classical approaches for problems with only continuous variables.  The latter two were proposed by the authors of the current paper and tested only on problems with continuous variables.  However, the weighted-constraint scalarization has been shown to be efficient for approximating Pareto points when the domain and/or front might be disconnected \cite{BurKayRiz2014, ParZilZil2017}.  Even though these techniques exist in the literature, for the convenience of the reader, we give here short descriptions. For more details, the reader is referred to \cite{BurKayRiz2014,BurKayRiz2017}, as well as \cite[Sections 2.4 and 3.2.4]{ParZilZil2017}, and the references therein.

\noindent
{\bf Weighted sum scalarization.}\ \
This type of scalarization was introduced by Gass and Saaty in \cite{Gass1955} and Zadeh in \cite{Zadeh1963}, it is computationally cheap and easy to implement. However, it cannot generate efficient points lying in nonconvex\footnote{The word {\em nonconvex} is loosely used here to describe those sections of the front which are ``dented", or, ``going inwards".} sections of the front.  The idea behind this method is to minimize a weighted sum  of all objectives. Namely, with fixed $w\in W$, the weighted-sum scalarization of Problem~$(P)$ is given by
\begin{equation} \tag{P$_w$}
\ds\min_{x \in X} \ \ \ds\sum_{i=1}^{\ell}w_i\,f_i(x)\,.
\end{equation}
Every solution of (P$_w$) is weak efficient (see
\cite{Miettinen1999}), and this fact is used in constructing an approximation of the Pareto front. For a fixed $w\in W$, define the solution set of P$_w$ as
\[
Sol(P_w) := \{ x \in X \mid x\ \mbox{solves}\ (P_w) \}\,.
\]

\noindent
{\bf The {\boldmath$k$}th-objective {\boldmath$\epsilon$}-constraint scalarization.}\ \
The method was introduced by Changkong and Haimes, and a comprehensive analysis can be found in \cite{ChaHai1983}. Let $\epsilon := (\epsilon_1,\ldots,\epsilon_{\ell})\in \mathbb{R}^{\ell}$ be fixed. In this approach, each of the objectives (say, the function $f_k$), is minimized and the remaining $(\ell-1)$ objectives are constrained by upper bounds $\epsilon_i \in \mathbb{R}$, $i=1,...,\ell$, $i \neq k$. The $k$th-objective $\epsilon$-constraint scalarization for Problem~$(P)$ is defined as
\begin{equation} \tag{$P^k(\epsilon)$}
 \left\{\begin{array}{cl} \ds \min_{x\in  X} & \ \ f_{k}(x) \\[2mm]
\mbox{s.t.} & \ \ f_i(x) \leq \epsilon_i\,,\ \ \forall\;i=1,\ldots,\ell\,, \;\; i \neq k,
\end{array}
\right.
\end{equation}
where $k \in \{1,2,...,\ell\}.$  Unlike the weighted-sum approach, this method can generate the Pareto points located in nonconvex sections of the front. However, the upper bounds $\epsilon_i$ in the constraints have to be chosen carefully, otherwise the new feasible region might be empty. In order to avoid this situation, a suitable range of values for $\epsilon_i$ has to	be known beforehand.  In \cite[p. 85--86]{Miettinen1999}, it is shown that if $\bar{x}$ solves ($P^k(\epsilon)$) for some $k$ then $\bar{x} \in WE(P)$.  For $\epsilon\in\mathbb{R}^{\ell}$, we define
\[
Sol(P^k(\epsilon)) := \{ x \in X \mid x\ \mbox{solves}\ (P^k(\epsilon)) \}\,.
\]
\noindent
{\bf Pascoletti--Serafini scalarization.}\ \ Introduced by Pascoletti and Serafini in~\cite{PasSer1984}, this approach is also referred to as the \emph{goal-attainment method} (see \cite{ColSia2004, DutKay2011, Miettinen1999, PasSer1984}).  Given a fixed $w \in W$, the Pascoletti--Serafini scalarization is posed as the following problem.
\begin{equation} \tag{$PS$}  \label{P-S}
 \left\{\begin{array}{cl} \ds \min_{(\alpha, x) \in \mathbb{R} \times X} &\ \ \alpha \\[2mm]
\mbox{s.t.} & \ \ w_i\,(f_i(x)-u_i) \leq  \alpha\,,\ \ \forall\;i=1,\ldots,\ell\,,
\end{array}
\right.
\end{equation}
where $\alpha \in \mathbb{R}$ is a new variable and $u=(u_1,\ldots, u_{\ell})$ is a utopia vector.  It should be noted that Problem~$(PS)$ can be reformulated as the (weighted) Tchebychev(-norm) scalarization~\cite{DutKay2011}. Every solution of Problem~$(PS)$ is a weak efficient point \cite{Eichfelder2008}. Hence, this technique is used for approximating connected Pareto fronts in~\cite{MueGraeSch2009}.  A modified version of the Pascoletti--Serafini scalarization has been used in~\cite{DasDen1998}, again for generating connected Pareto fronts. \\

\noindent
{\bf Weighted-constraint scalarization.}\ \ Introduced in~\cite{BurKayRiz2014} by the authors of the current paper, this scalarization has been numerically illustrated to be useful in finding weak efficient points in a disconnected Pareto front, in the presence of a disconnected feasible set or a disconnected domain. As in the $k$th objective $\epsilon$-constraint approach, each objective $f_k$, $k=1,\ldots,\ell$, is minimized separately; however, all of the objectives appear in the constraints of the scalarized problem.  Namely, the following $\ell$ problems are solved.
\begin{equation} \tag{$P_w^k$}
\left\{\begin{array}{cl} \ds \min_{x\in  X} & \ \ f_k(x) \\[2mm]
\mbox{s.t.} & \ \ w_i\,f_i(x) \leq w_k\,f_k(x)\,,\quad\forall\;i=1,\ldots,\ell, \;\; i \neq k\,,
\end{array}
\right.
\end{equation}
where $k = 1,\ldots,\ell$.  Each Problem $(P_w^k)$, for $k = 1,\ldots,\ell$, is referred to as a {\em subproblem} of the scalarization. For fixed $w\in W$ and $k\in\{1,\ldots,\ell\}$, define the solution set of a subproblem as
\[
Sol(P_w^k) :=  \left \{x \in X \mid x \;\; \mbox{solves}\;\; (P_w^k) \right \}\,.
\]
Given $w'\in W$, it is proved in  \cite{BurKayRiz2014} that 
\begin{equation}\label{pkw}
\ds \bigcap_{k=1}^l Sol(P_{w'}^k)\subset WE(P)\subset \bigcup_{w\in W} \left[\ds \bigcap_{k=1}^l Sol(P_w^k)\right]\,.
\end{equation}
When the leftmost expression is nonempty for some $w'\in W$, the leftmost inclusion above can be used to obtain points in the Pareto front. If such a weight $w'$ is not directly available then \cite[Theorem 3.1 and Proposition 3.3]{BurKayRiz2014} can be used to generate new points in the front. To make these statements precise, we recall next the relevant results. The first of these results ensures that the map from the space of weights $W$ to the Pareto set $WE(P)$ is surjective.  This property allows one to construct an algorithm using some discretization/partition of $W$, which can be utilized to obtain an approximation of the Pareto front. 

\begin{theorem}(\cite[Theorem~3.1]{BurKayRiz2014}) \label{surjectivity}
A point $x \in X$ is a weak efficient solution of Problem~$(P)$, if and only if there exists some $w \in W$ such that $x\in Sol(P_w^k)$, for all $k \in\{1,\ldots,\ell\}$. If $x\in WE(P)$, the required $w$ is given by
\[
w_i := \ds\frac{1/f_i({x})}{\ds\sum_{j=1}^{\ell}1/f_j({x})}\,,\quad i=1,\ldots,\ell\,.
\]
\end{theorem}
Theorem~\ref{surjectivity} implies that, for a fixed $w\in W$, if
\[
x\in \ds\bigcap_{k=1}^\ell Sol(P_w^k)\,,
\]
then $x$ is a weak efficient point. On the other hand, if solutions of two subproblems are different from one another then a comparison is made between these solutions to eliminate a dominated point. If a new solution of a subproblem is non-dominated compared with the solutions of the other subproblems,  then that solution is also weak efficient.  Proposition~\ref{different_notcomparable} below states this fact more precisely and plays an important role in the weighted-constraint scalarization approach. The idea in the proposition is implemented in Algorithms~5--7 in the current article for the removal of dominated points.

\begin{proposition}(\cite[Proposition~3.3]{BurKayRiz2014}) \label{different_notcomparable}
\ Assume that $\exists\, w \in W$ such that $Sol(P_w^j) \neq \emptyset$,\linebreak $\forall j=1,\ldots,\ell $.  Suppose that, for some $k \in\{1,\ldots,\ell\}$, $\exists\, \bar{x}_k \in Sol(P_w^k)$ such that $\forall r \in \{1,2, \ldots, \ell \}$, $r \neq k$, $\exists\,\bar{x}_r \in Sol(P_w^r)$, which satisfies
\begin{equation}\label{corr1}
  f_r(\bar{x}_r) \geq  f_r(\bar{x}_k)\,.
\end{equation}
Then $\bar{x}_k \in WE(P)$.
\end{proposition}

\noindent{\bf The feasible-value-constraint scalarization.}\ \
Introduced in \cite{BurKayRiz2017} by the authors of the current paper, this scalarization technique has numerically been illustrated to efficiently approximate the boundary and the interior of the Pareto front even when the front is disconnected. It uses the specific expression of the weights given in Theorem~\ref{surjectivity}, evaluated at a feasible point. Namely that with some $\hat{x} \in X$, one sets
\begin{equation}\label{weiassum}
w_i:=\frac{1/f_i(\hat{x})}{\ds\sum_{j=1}^{\ell}1/f_j(\hat{x})}\,,\quad i=1,\ldots,\ell\,.
\end{equation}
For $w \in W$, we obtain
\begin{equation}\label{property}
w_kf_k(\hat{x})=w_if_i(\hat{x}),\ \ \forall i=1,2,..., \ell,\ \ i \neq k.
\end{equation}
We will usually assume that a weight satisfies \eqref{weiassum} (and hence \eqref{property} holds).  The associated scalar problem is defined as
\begin{equation} \tag{$P_{\hat{x}}^k$}
\ \left\{\begin{array}{cl} \ds\min_{x\in  X} & \ \ w_{k}f_{k}(x),
\\[4mm]
\mbox{s.t.} & \ \ w_if_i(x) \leq w_kf_k(\hat{x}),\ \ i=1,2,...,\ell\,,\ i \neq k.
\end{array}
\right.
\end{equation}
For fixed $w\in W$, $\hat{x}\in X$ and $k\in\{1,\ldots,\ell\}$, define the solution set
\[
Sol(P_{\hat{x}}^k) := \{ x \in X \mid x\ \mbox{solves}\ (P_{\hat{x}}^k) \}\,.
\]
If $\bar{x} \in Sol(P_{\hat{x}}^k)$ for some $k$, then $\bar{x} \in WE(P)$ \cite[Theorem 4.6]{BurKayRiz2017}. Note that this method is reminiscent of the $k$th-objective $\epsilon$-constraint approach, but uses a very specific choice of the vector $\epsilon$, the components of which are determined by the feasible objective function values. The relationship between these two approaches is made precise in the next section.

\subsection{The relationship between {\boldmath$(P^k(\epsilon))$} and {\boldmath$(P_{\hat{x}}^k)$}}
\label{sec3.1}
	
Given $\hat{x}\in X$, and $w$ as in \eqref{weiassum}, we show next that if $\epsilon\le \hat{\epsilon} := f(\hat{x})$, then the set of solutions of Problem~($P^k(\epsilon)$) contains the set of solutions of Problem~($P^k_{\hat{x}}$). Recall that \eqref{positive} yields $f(x)\ge 0$ for every $x\in X$, so $\hat\epsilon\in \mathbb{R}^{\ell}_+$.

\begin{theorem}\label{sc1_rela1}
Fix $\hat{x}\in X$, let $w$ be as in \eqref{weiassum}, and set $\hat\epsilon:=f(\hat x)$. 
Then $Sol(P_{\hat{x}}^k)={Sol}(P^k(\hat \epsilon))$. 

\end{theorem}
\begin{proof}
Note that the objective functions of both problems are equivalent because $w\in W$. Regarding the constraints, we can use \eqref{property} and the definition of $\hat\epsilon$ to write $w_kf_k(\hat{x})=w_if_i(\hat{x})=w_i \hat\epsilon_i$ for every $i \neq k$. Using again the fact that $w\in W$, the constraints in Problem~($P^k(\hat\epsilon)$) can be equivalently written as
\[
w_i f_i(x)\le w_i \hat\epsilon_i= w_i f_i(\hat x) = w_k f_k(\hat{x}),
\]
which are the constraints of Problem ($P_{\hat{x}}^k$). Hence, with this choice of $\hat\epsilon$, Problems~($P^k(\hat \epsilon)$) and ($P_{\hat{x}}^k$) are equivalent. In particular, they have the same solutions. Namely that ${Sol}(P^k(\hat \epsilon)) = Sol(P_{\hat{x}}^k)$\,.
\end{proof}

\subsection{The relationship between {\boldmath$(P_w)$} and {\boldmath$(P_{\hat{x}}^k)$}}
\label{sec3.2}

In Theorem~\ref{sc1_rela3} below, we show that those Pareto points which can be generated by the weighted-sum scalarization can also be generated by the feasible-value-constraint scalarization. The converse of this statement, however, is not true: those Pareto points which can be generated by the feasible-value-constraint approach may not be attainable by the weighted-sum method. We recall that every Pareto point is a solution of Problem~$(P_{\hat{x}}^k)$ by Theorem~\ref{surjectivity}, but Problem~$(P_w)$ cannot generate some of the Pareto points located in a part/region of the front whose epigraph is a nonconvex set.
\begin{theorem}\label{sc1_rela3}
Fix $w\in W$ and $\bar{x}\in Sol(P_w)$.  Then $\bar{x}\in Sol(P_{\bar{x}}^j)$ holds for all $j$.
\end{theorem}
\begin{proof} 
Let $\alpha \in W$ be as in \eqref{weiassum}-\eqref{property} for $\hat{x} = \bar{x}$. Since $\bar{x} \in Sol(P_w)$, we have that
\begin{equation}\label{sc1_equ3}
\sum_{i=1}^{\ell}w_i[f_i(x)-f_i(\bar{x})] \geq 0, \;\;
\forall x \in X.
\end{equation}
Let us note first that $\bar{x}$ is feasible for $(P^j_{\bar{x}})$ because the constraint set is defined using the weight vector $\alpha$. Assume that, for some $j \in\{1,\ldots,\ell\}$,  $\bar{x} \notin Sol(P_{\bar{x}}^j)$ . Then there exists $x \in X$ such that
\begin{equation} \label{rel_pw_pk1}
\alpha_jf_j(x)<\alpha_jf_j(\bar{x}),
\end{equation}
and
\begin{equation}\label{rel_pw_pk2}
\alpha_if_i(x)\leq \alpha_jf_j(\bar{x})=\alpha_if_i(\bar{x}), \;\; i=1,2,...,\ell \;\;\; \mbox{and} \;\;\;i \neq j,
\end{equation}
where we have used the definition of $\alpha$ in the equality. From \eqref{rel_pw_pk1} and \eqref{rel_pw_pk2}, we have that
\begin{equation} \label{rel_pw_pk3}
\hspace*{20mm} f_j(x)-f_j(\bar{x})<0,
\end{equation}	
and 
\begin{equation} \label{rel_pw_pk4}
f_i(x)-f_i(\bar{x}) \leq 0, \;\; i\neq j\,,
\end{equation}	
 where $x$ and $j$ are as in \eqref{rel_pw_pk1}.

Since $ w _i >0 $, and using \eqref{rel_pw_pk3} and \eqref{rel_pw_pk4}, we conclude that
	\[
	\hspace*{20mm} \sum_{i =1}^{\ell}w_i(f_i(x)-f_i(\bar{x}))<0,
	\]
	which contradicts \eqref{sc1_equ3}. Therefore $\bar{x} \in
	Sol(P_{\bar{x}}^j)$ for all $j$.
\end{proof}

\begin{remark}\rm
If we assume that $w\ge0$, instead of $w>0$, a reasoning similar to the one in the proof above cannot be used directly in getting the conclusion of Theorem~\ref{sc1_rela3}. However, if we assume further that $Sol(P_w)$ is a singleton for some weight $w\ge 0$, a trivial modification of the proof will work.  Since the modified proof is straightforward, we do not elaborate it here.
\end{remark}



\subsection{The relationship between {\boldmath$(P^k(\epsilon))$} and {\boldmath$(P_w^k)$}}
\label{sec3.3}

Theorem~\ref{surjectivity} states that every weak efficient point is the solution of the $k$th-objective weighted-constraint scalarization problem~$(P_w^k)$, for a suitable choice of weight. On the other hand, some weak efficient points may not be  solutions of the $k$th-objective $\epsilon$-constraint problem~$(P^k(\epsilon))$, if $\epsilon$ is too large. 

\begin{theorem}\label{sc2_epsilon}
	Assume that $\bar{x}, \hat{x} \in WE(P)$ and  
	\begin{equation}\label{sc1_eps1}
	f_k(\hat{x})< f_k(\bar{x})\;\; \mbox{for some}\; k.
	\end{equation}
 If $\epsilon \geq \hat{\epsilon}:=f(\hat{x})$, then  $\bar{x} \notin Sol(P^k(\epsilon)).$ 
\end{theorem}
\begin{proof} Assume that there exists $\tilde{\epsilon} \geq \hat{\epsilon}$ such that $\bar{x} \in Sol(P^k(\tilde{\epsilon}))$. Thus, we can write
\begin{equation}\label{sc1_eps2}
f_k(\bar{x})\leq f_k(x),\;\; \forall\;x\in X \; \mbox{such that}\; f_i(x)\leq \tilde{\epsilon}_i, \;  i\neq k.
\end{equation}	
Note that $f_i(\hat{x})=\hat{\epsilon}_i\le \tilde{\epsilon}_i$ for all $i \neq k$. Hence, $x=\hat{x}$ satisfies the constraints of problem $P^k(\tilde{\epsilon})$. Therefore, using $x=\hat{x}$ in \eqref{sc1_eps2}, gives 
	\[
	f_k(\bar{x})\leq f_k(\hat{x}),
	\]
contradicting \eqref{sc1_eps1}. Therefore, for all $\epsilon \geq \hat{\epsilon}$,  we must have $\bar{x} \notin Sol(P^k(\epsilon))$. 
\end{proof}

\begin{remark} \label{rem:choice} \rm
Recall that Theorem~\ref{sc1_rela1} implies, with the choice of $\epsilon = f_k(\hat{x})$ where $\hat{x}\in X$, that the solution sets of $(P^k(\epsilon))$ and $(P_{\hat{x}}^k)$ are the same.  On the other hand, Theorem~\ref{sc2_epsilon} asserts that, if the Pareto front has points which are weak efficient (but not just efficient) then, unless $\epsilon$ is chosen in an informed manner, the solution set of $(P^k(\epsilon))$ may not contain all of these weak efficient points.  However, by Theorem~\ref{surjectivity}~\cite[Theorem~3.1]{BurKayRiz2014}), the solution set of $(P_w^k)$ contains the whole weak efficient set.  This justifies the employment of the weighted-constraint scalarization~$(P_w^k)$ in the algorithms in the next section, rather than the scalarizations~$(P^k(\epsilon))$, $(P_{\hat{x}}^k)$, or indeed the weighted sum scalarization~$(P_w)$.
\end{remark}

\section{Algorithms}
\label{sec4}

\subsection{Grid generation techniques}
In our study, we incorporate two types of grid generation techniques in designing the algorithms for approximating Pareto points. For  details of these grid generation processes, the reader is referred to \cite{DasDen1998,MueGraeSch2009}, as well as \cite[Sections~1.2 and 5]{BurKayRiz2017}.

\noindent
{\bf Convex Hull of the Individual Minima (CHIM).}\ \ The CHIM grid generation technique was proposed by Das and Dennis, and employed as part of their Normal Boundary Intersection (NBI) method in~\cite{DasDen1998}.  The NBI method is arguably the most popular approach to constructing an approximation of the Pareto front. First, the vector of individual minima is obtained, and then the convex hull of these individual minima is generated. After a {\em uniform} discretization of the convex hull, and assigning a weight to each node of the discretization, the CHIM grid is created.

In the NBI method, given the CHIM grid along with the corresponding array of weights, a modified Pascoletti--Serafini scalarization of the problem is solved for each weight in the array to generate an approximation of the Pareto point. Geometric illustrations of the CHIM grid generation process for two and three objectives can be found in~\cite[Figure 1]{BurKayRiz2017}.

\noindent
{\bf Sequential Boundary Generation (SBG).}\ \ The SBG grid generation technique was proposed by Mueller-Gritschneder et al.\ in~\cite{MueGraeSch2009}. We observe in our experiments that the SBG method works well when the CHIM grid is unable to generate the true boundary of the Pareto front. The SBG method constructs the Pareto front of the problem sequentially. In doing this, the SBG method solves linear programming problems. As a result, it requires more computational time than the CHIM  method. A detailed description of the SBG technique can be found in \cite[Sections~1.2 and 5]{BurKayRiz2017}.

As exemplified via test problems in \cite{BurKayRiz2017}, when the boundary of the Pareto front happens to be complicated (which is a common occasion with the three- and four-objective problems), the SBG grid is more suitable than the CHIM grid in getting a correct depiction of the boundary.

\subsection{Descriptions of the algorithms}
\label{descriptions}

Table~\ref{algorithms} provides a list of the algorithms we implement to solve two-, three- and four-objective problems.  In the table, we list the number of objectives that an algorithm can handle, as well as the grid generation and the scalarization techniques that are utilized in that algorithm.  Algorithms~1 and 2 can solve two-objective problems by using the CHIM grid.  While Algorithms~3-6 can all solve three-objective problems, Algorithms~3 and 4 employ the CHIM grid generation technique, and Algorithms~5 and 6 employ the SBG grid generation technique.  Algorithm~7 can handle four objectives and utilizes the SBG grid generation  technique.  As for the scalarization techniques, while the odd-numbered algorithms use the weighted-constraint scalarization, as in Problem~(P$_w^k$), the even numbered ones use the Pascoletti--Serafini scalarization, as in Problem~(PS).  

As indicated in the last column of Table~\ref{algorithms}, Algorithms 1 and 2 directly come from \cite{BurKayRiz2014}, and Algorithms 3--6 from \cite{BurKayRiz2017}.  It should be noted however that Algorithm~3 in the current paper was earlier provided in \cite{BurKayRiz2017} as Algorithm~2 with steps written in an abridged form. Here we provide a fully expanded and explicit version of that algorithm as Algorithm~3 in Appendix~A.1.   Algorithm~7, on the other hand, does not appear elsewhere, and we provide its full description in Appendix~A.2.  We will use Algorithm~7 particularly for solving the challenging rocket injector design problem in Section~\ref{rocket}.

\begin{table}[t]
  \caption{\small{\textit{The algorithms implemented in the paper.}}}
\footnotesize
\vskip 1.5em
\centering
\begin{tabular}{|c| c| c| c| c|}
\hline
 & Number of & Grid & &  \\
Algorithm & objectives & generation & Scalarization & Reference \\ \hline
1 & 2 & CHIM & Weighted-constraint & \cite{BurKayRiz2014} \\[1mm]
2 & 2 & CHIM & Pascoletti--Serafini & \cite{BurKayRiz2014} \\[1mm] \hline
3 & 3 & CHIM & Weighted-constraint & \cite{BurKayRiz2017} \\[1mm]
4 & 3 & CHIM & Pascoletti--Serafini & \cite{BurKayRiz2017} \\[1mm]
5 & 3 & SBG & Weighted-constraint & \cite{BurKayRiz2017} \\[1mm]
6 & 3 & SBG & Pascoletti--Serafini & \ \ \ \ \cite{BurKayRiz2017, MueGraeSch2009} \\[1mm] \hline
7 & 4 & SBG & Weighted-constraint & New \\[1mm]
\hline
\end{tabular}
\label{algorithms}
\end{table}

\section{Numerical Experiments}
\label{sec5}

In this section, we test and compare Algorithms~1--7 that are described in Section~\ref{sec4}, by means of two-, three- and four-objective integer and mixed-integer programming problems. The test problems~1--3 are integer programming problems.  These problems are designed in such a way that the number of points in the Pareto front is finite and that they can be interpreted (or visualized easily) geometrically, so that we know the weak Pareto points prior to computations. These known weak Pareto, or weak efficient, points are referred to as being {\em exact}.  One of our aims is to understand the capabilities of Algorithms~1--6 in approximating the set of exact weak Pareto points. The task of approximating a Pareto front is particularly challenging for the three- and four-objective cases.  We test Algorithm~7 on a challenging real-life problem, namely the rocket injector design, which has four objective functions to minimize simultaneously.

We have coded all of Algorithms~1--7 in MATLAB, and utilized the optimization software BARON \cite{Sahinidis2017, TawSah2005} or SCIP \cite{SCIP2017, Achterberg2009} or Bonmin \cite{Bonami2008} or IPOPT \cite{Wachter2006}, for solving the single-objective integer and mixed-integer subproblems (P$_w^k$), including the scalarized problems, in each algorithm.  In particular, we used BARON in solving all the test problems. Additionally, we used SCIP, Bonmin and IPOPT for the last test problem, rocket injector design, in order to reconfirm the numerical results we obtained. All solver tolerances are set at $10^{-10}$. In the generation of the SBG grid, MATLAB's {\tt linprog} was used, with default options, for solving the associated linear programming~(LP) problems. The computations have been performed on an HP ENVY 14 laptop with 4~GB RAM and core i5 at 1.6 GHz.

\subsection{Test problem 1: Two objective functions} 
\label{testproblem1}

\begin{figure}[t!]
\begin{minipage}{80mm}
\begin{center}
\hspace*{-5mm}
\psfrag{f1}{$f_1$}
\psfrag{f2}{$f_2$}
\includegraphics[width=80mm]{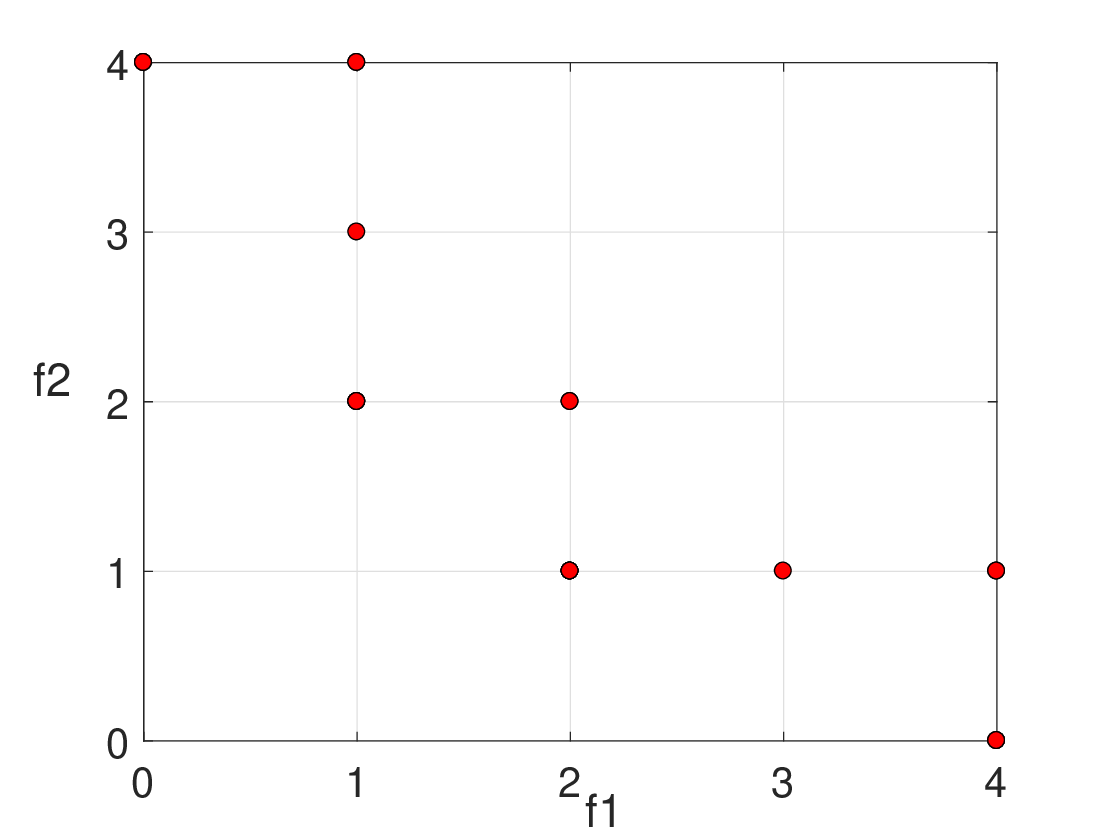} \\[3mm]
{\footnotesize (a) Pareto front by Algorithm~1}.\\
\end{center}
\end{minipage}
\begin{minipage}{80mm}
\begin{center}
\hspace*{-5mm}
\psfrag{f1}{$f_1$}
\psfrag{f2}{$f_2$}
\includegraphics[width=80mm]{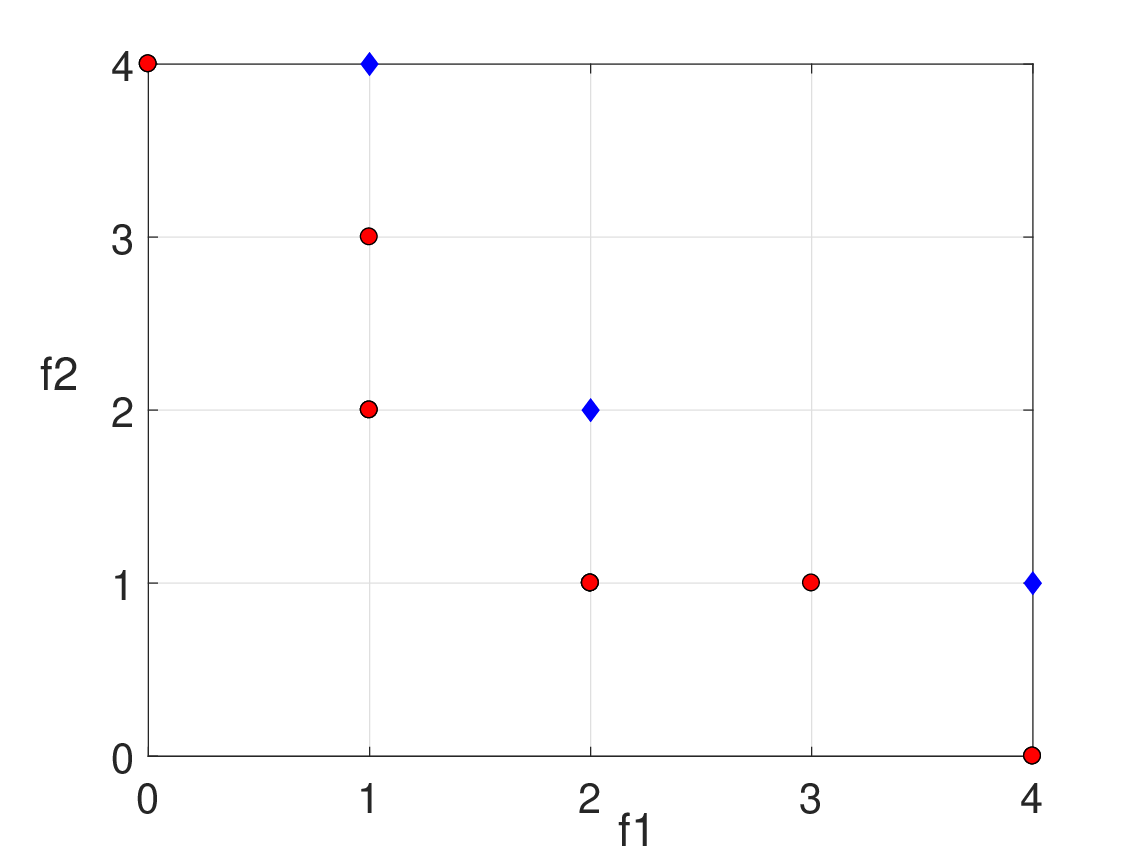} \\[3mm]
{\footnotesize (b) Pareto front by Algorithm~2}.
\end{center}
\end{minipage}
\caption{\sf\small Test problem 1 -- Pareto fronts constructed by Algorithms~1 and 2, shown by (red) circles. The (blue) diamonds represent the weak Pareto points that were missed by Algorithm~2.}
\label{figexa_1}
\end{figure}

We consider, mainly for illustration purposes, a small problem with two objective functions and two integer variables:
\[
\begin{array}{cl} \rm\min & \ (x_1,\ x_2) \\[2mm]
\mbox{s.t.} & (x_1-4)^2 + (x_2-4)^2 -16 \leq 0\,, \\[1mm]
  & 0\leq x_1,x_2\leq 4\,,\;\; \mbox{and}\;\;x_1, x_2 \mbox{ are integers}.
\end{array}
\]
The problem above is easy to interpret: The feasible points of the above problem are those points with integer coordinates in the lower-left quarter of a circle of radius 4 centred at $(4,4)$.  One can easily count that there are 17 feasible points of the problem.  Then, of these feasible points, the ones which are located as far to the ``south-west" as possible are weak Pareto points.  There are nine weak Pareto points, the set of which can simply be written as (see Figure~\ref{figexa_1}(a))
\[
\{(0,4),(1,2),(1,3),(1,4),(2,1),(2,2),(3,1),(4,0),(4,1)\}.
\]
The above set is nothing but the Pareto front, which is discrete and so is disconnected.  Note that $f_1(x)=x_1$ and $f_2(x)=x_2$ are the objective function of the problem.  As the utopia point, we take $u = (-10, -10)$.  For Algorithms~1 and 2, we use the CHIM grid. Since the CPU time requirements are in general different for each of the algorithms for the same number of the grid points, we have adjusted the number of grid points in such a way that the CPU time each algorithm takes is roughly the same (in this case, around 13 and 14 seconds, respectively).  This way, we provide the two algorithms an equal footing in the comparisons. The data on the CPU time, the number of subproblems attempted, and the number of weak Pareto points generated by each algorithm, are reported in Table~\ref{table:exp-1a}.  Parts~(a) and (b) of Figure~\ref{figexa_1} display all of the Pareto points, i.e., an approximation of the Pareto front, that could be obtained by each algorithm, respectively, within the CPU times indicated.

Algorithm 1, as shown by Figure~\ref{figexa_1}(a), successfully finds all nine weak Pareto points in the Pareto front. On the other hand, Algorithm 2 is able to generate only 6 of these points. When we increase the number of grid points, i.e., increase the allowed CPU time, by (say) ten times, the number of Pareto points obtained by Algorithm~2 still turns out to be 6.  In other words, no improvement in the approximation of the Pareto front can be obtained by Algorithm~2, even after providing the algorithm with a much finer grid and so allowing a much longer computational time.

\begin{table}[t]
  \caption{\small{\textit{Test problem 1 -- Numerical performance of Algorithms~1 and 2.}}}
\footnotesize
\vskip 1.5em
\centering
\begin{tabular}{|c| c| c| c| c|}
\hline
		&  			& Number of 	& Number of  \\
                 & CPU time	& subproblems & weak Pareto  \\
 Algorithm &     [sec] 		&  attempted   	& points generated   \\ [0.5ex]
\hline
1 &14 & $11 \times 2 = 22$ & 9 (all)  \\ [.5ex]\hline
2 &13& 14 & 6\ \ \ \ \ \ \  \\ [.5ex]
\hline
\end{tabular}
\label{table:exp-1a}
\end{table}

\subsection{Test problem 2: Three objective functions -- nonlinear constraints}
\label{tri-obj}

Next we modify the test problem~1 in Section~\ref{testproblem1} by introducing one more objective function and one more variable as follows.
\[
\begin{array}{cl} \min & \ \left[x_1, x_2, x_3\right] \\[2mm]
 \mbox{s.t} &\  (x_1-2)^2+(x_2-2)^2+(x_3-2)^2 \leq 4, \\[1mm]
    &  \ 0 \leq x_i \leq 4\,,\ i=1,2,3,\mbox{ are integers.}
\end{array}
\]

\begin{figure}[t!]
\hspace{-1cm}
\begin{minipage}{90mm}
\begin{center}
\hspace*{0cm}
\psfrag{f1}{$f_1$}
\psfrag{f2}{$f_2$}
\psfrag{f3}{$f_3$}
\includegraphics[width=80mm]{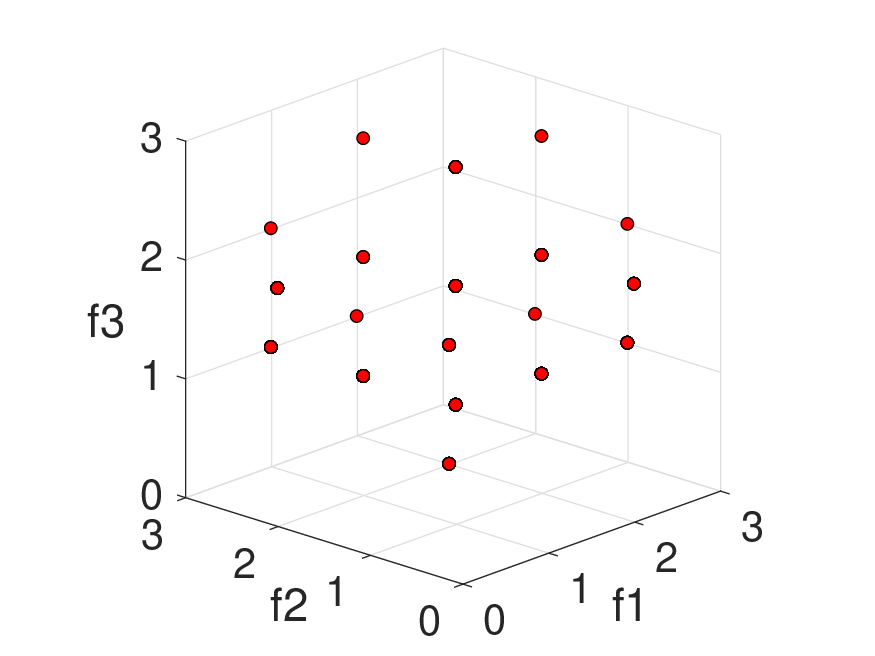} \\
{\footnotesize (a) Pareto front by Algorithm 3.}
\end{center}
\end{minipage}
\hspace{-1cm}
\begin{minipage}{85mm}
\begin{center}
\hspace*{0cm}
\psfrag{f1}{$f_1$}
\psfrag{f2}{$f_2$}
\psfrag{f3}{$f_3$}
\includegraphics[width=80mm]{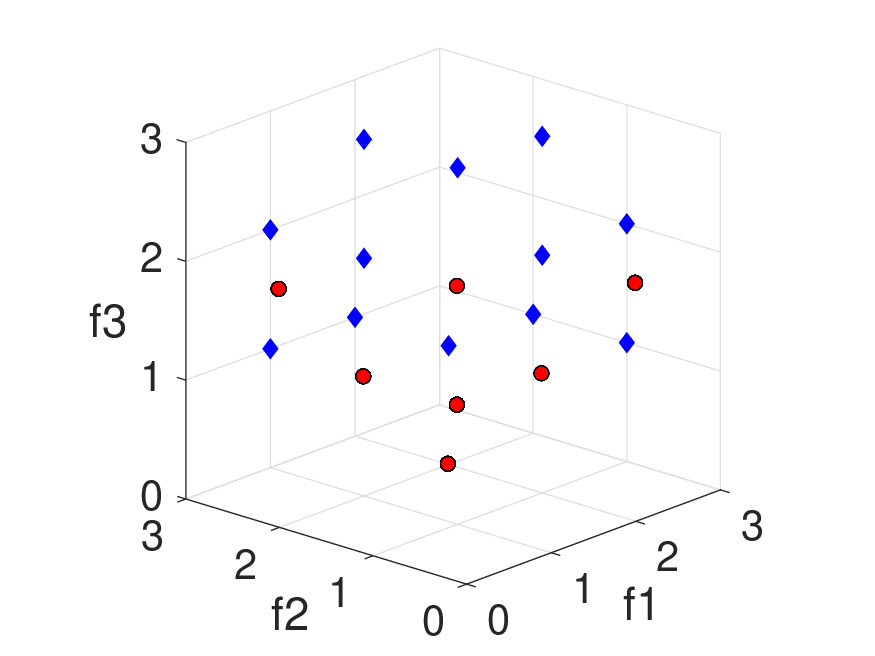} \\
{\footnotesize (b) Pareto front by Algorithm 4.}
\end{center}
\end{minipage}
\\[2mm]
\hspace*{-1cm}
\begin{minipage}{90mm}
\begin{center}
\hspace*{0cm}
\psfrag{f1}{$f_1$}
\psfrag{f2}{$f_2$}
\psfrag{f3}{$f_3$}
\psfrag{(1,2,3)}{\scriptsize $(1,2,3)$}
\psfrag{(1,1,1)}{\scriptsize $(1,1,1)$}
\psfrag{(1,2,1)}{\scriptsize $(1,2,1)$}
\psfrag{(1,3,1)}{\scriptsize $(1,3,1)$}
\psfrag{(2,2,0)}{\scriptsize $(2,2,0)$}
\psfrag{(2,0,2)}{\scriptsize $(2,0,2)$}
\includegraphics[width=80mm]{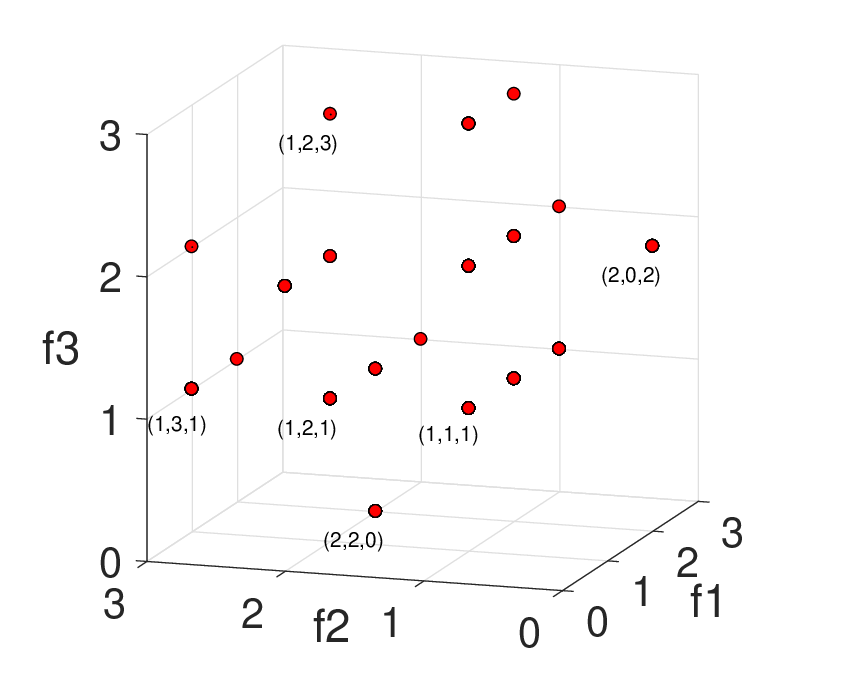} \\
{\footnotesize (c) A rotated view of the front in (a).}
\end{center}
\end{minipage}
\hspace{-1cm}
\begin{minipage}{85mm}
\begin{center}
\hspace*{0cm}
\psfrag{f1}{$f_1$}
\psfrag{f2}{$f_2$}
\psfrag{f3}{$f_3$}
\psfrag{(1,2,3)}{\scriptsize $(1,2,3)$}
\psfrag{(1,1,1)}{\scriptsize $(1,1,1)$}
\psfrag{(1,2,1)}{\scriptsize $(1,2,1)$}
\psfrag{(1,3,1)}{\scriptsize $(1,3,1)$}
\psfrag{(2,2,0)}{\scriptsize $(2,2,0)$}
\psfrag{(2,0,2)}{\scriptsize $(2,0,2)$}
\includegraphics[width=80mm]{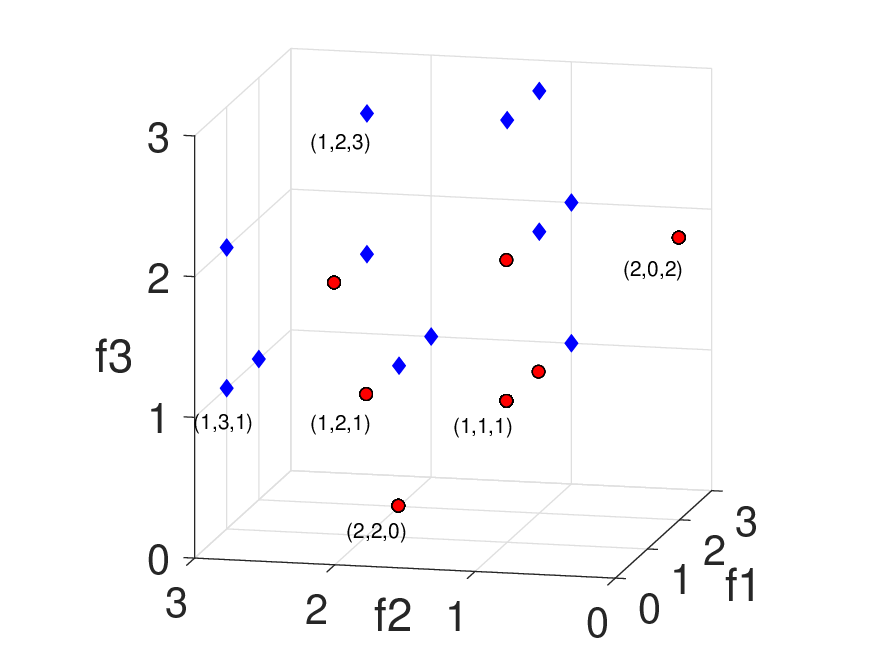} \\
{\footnotesize (d) A rotated view of the front in (b).}
\end{center}
\end{minipage}
\caption{\sf\small Test problem 2 -- Pareto fronts constructed by Algorithms~3 and 4, shown by (red) circles. The (blue) diamonds represent the weak Pareto points that were missed by Algorithm~4.}
\label{testpro2I}
\end{figure}

The feasible region is now defined to be a spherical ball rather than a part of a circular region; namely, the feasible set consists of the points with integer coordinates in a sphere of radius 2 centred at $(2,2,2)$.  In this case, the problem has 70 feasible points, of which 19 are weak Pareto.  The set of weak Pareto points, or the Pareto front, can simply be written as  (see Figures~\ref{testpro2I}(a) and \ref{testpro2I}(c))
\begin{eqnarray*}
&& \{(0,2,2),(1,1,1),(1,1,2),(1,1,3),(1,2,1),(1,2,2),(1,2,3),(1,3,1),(1,3,2), \\[1mm]
&&\ \,(2,0,2),(2,1,1),(2,1,2),(2,1,3),(2,2,0),(2,2,1),(2,3,1),(3,1,1),(3,1,2),(3,2,1)\}.
\end{eqnarray*}

We note that $f_1(x)=x_1$, $f_2(x)=x_2$ and $f_3(x)=x_3$.  We choose the utopia point as $u = (-10,-10,-10)$.  In order to solve the problem, i.e., to obtain an approximation of the Pareto front, we have implemented Algorithms~3 and 4, both of which use the CHIM grid.  As in the case of Test Problem 1, with the choice of the same CHIM grid, each algorithm requires a substantially different length of CPU time.  Therefore, we have adjusted the number of points in the CHIM grids of the algorithms in such a way that the CPU time each algorithm takes is roughly the same (in this case, about 150 seconds each), so that the two algorithms can be compared on an equal footing.  Table~\ref{table:exp-2a} summarizes the performance of each algorithm.

\begin{table}[t!]
\caption{\small{\textit{Test problem 2 -- Numerical performance of Algorithms~3 and 4.}}}
\footnotesize
\vskip 1.5em
\centering
\begin{tabular}{|c| c| c| c| c|}
\hline
		&  			& Number of 	& Number of  \\
                 & CPU time	& subproblems & weak Pareto  \\
 Algorithm &     [sec] 		&  attempted   	& points generated   \\ [0.5ex]
\hline
3  		&147 		& $44 \times 3 = 132$ & 19 (all)  \\ [.5ex]\hline
4 		&150 		& 135 		& 7\ \ \ \ \ \  \\ [.5ex]
\hline
\end{tabular}
\label{table:exp-2a}
\end{table}

As can be seen in Table~\ref{table:exp-2a}, Algorithm~4 is able to find only 7 weak Pareto points out of 19, as can also be seen in Figures~\ref{testpro2I}(b) and (d).  Note that the Pareto front in Figure~\ref{testpro2I}(d) is a rotated view of the same front in Figure~\ref{testpro2I}(b) and has been added to the figure to make the visualization easier. Next we provided Algorithm~4 with much finer CHIM grids, which obviously resulted in far longer computational times; however, not a single new weak Pareto point could be obtained by Algorithm~4.  This indicates that Algorithm~4 is not efficient for this test problem. On the other hand, Algorithm~3 generates all the 19 weak Pareto points (see Figures~\ref{testpro2I}(a) and \ref{testpro2I}(c)) in 147 seconds of CPU time.

\subsection{Test problem 3: Three objective functions -- linear constraints}
\label{testproblem3}

This time, we consider the test problem given by Antunes et al.~\cite[Exercise~5, Section~6.5]{AntAlvCli2016}, who studied certain properties of two known weak Pareto points, but not the whole Pareto front. We are interested here in constructing the whole Pareto front of this problem, which, to the best of our knowledge, has not been done elsewhere, yet.
\[
\begin{array}{cl} \min & \ \left[-x_1, -x_2, -x_3\right] \\[2mm]
 \mbox{s.t} &\  3x_1+2x_2+3x_3 \leq 18, \\[1mm]
 &\  x_1+2x_2+x_3 \leq 10, \\[1mm]
 &\  9x_1+20x_2+7x_3 \leq 96, \\[1mm]
 &\  7x_1+20x_2+9x_3 \leq 96, \\[1mm]
 &  \  x_i \geq 0\,,\ i=1,2,3,\mbox{ are integers}.
\end{array}
\]
The problem has three non-negative integer decision variables.  The feasible region is defined by the intersection of four closed half-spaces as indicated above and has 83 feasible points (with non-negative integer coordinates), of which 60 are weak Pareto. The (disconnected) Pareto front consists of these 60 discrete points -- see the graphs in parts (a) and (c) of Figure~\ref{testpro1a}. 

\begin{figure}[hbt!]
\hspace{-1cm}
\begin{minipage}{90mm}
\begin{center}
\hspace*{0cm}
\psfrag{f1}{$f_1$}
\psfrag{f2}{$f_2$}
\psfrag{f3}{\hspace*{-1mm}$f_3$}
\includegraphics[width=75mm]{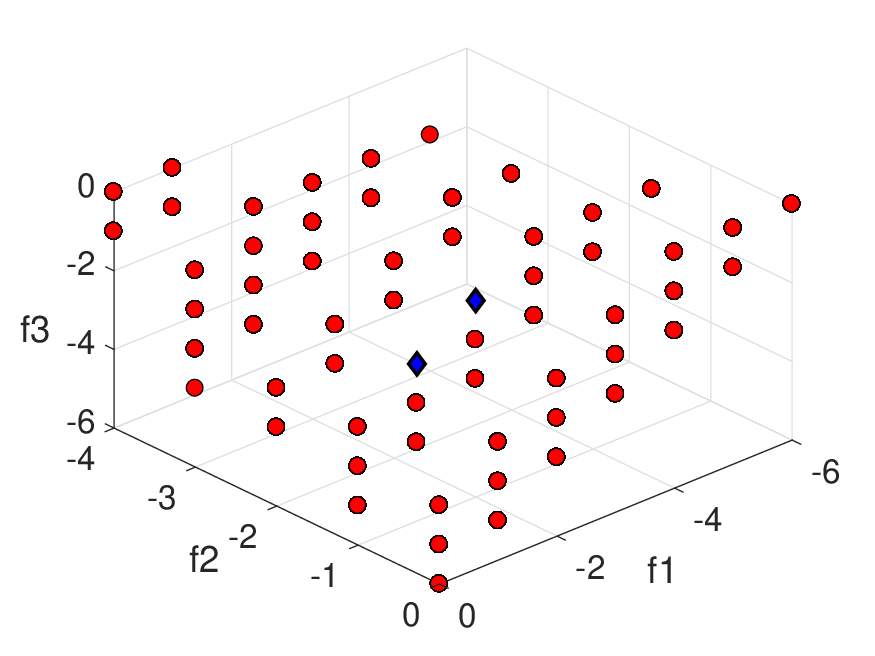} \\
{\footnotesize (a) Pareto front by Algorithm 3.}
\end{center}
\end{minipage}
\hspace{-1cm}
\begin{minipage}{90mm}
\begin{center}
\hspace*{0cm}
\psfrag{f1}{$f_1$}
\psfrag{f2}{$f_2$}
\psfrag{f3}{\hspace*{-1mm}$f_3$}
\includegraphics[width=75mm]{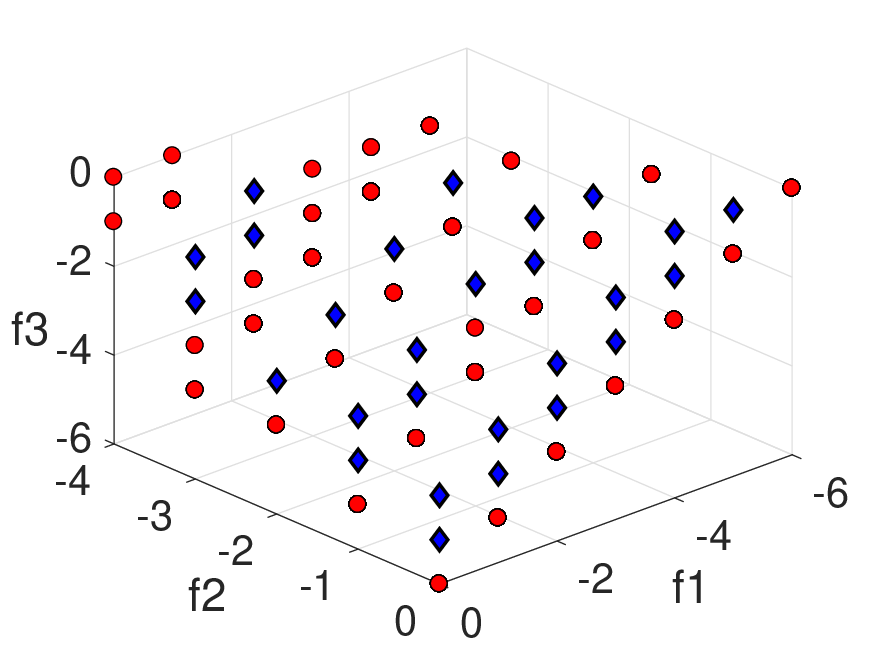} \\
{\footnotesize (b) Pareto front by Algorithm 4.}
\end{center}
\end{minipage}
\\[2mm]
\hspace*{-1cm}
\begin{minipage}{90mm}
\begin{center}
\hspace*{0cm}
\psfrag{f1}{$f_1$}
\psfrag{f2}{$f_2$}
\psfrag{f3}{\hspace*{-1mm}$f_3$}
\psfrag{(-1,-3,0)}{\tiny (-1,-3,0)}
\psfrag{(-3,-2,0)}{\tiny (-3,-2,0)}
\psfrag{(-4,0,0)}{\tiny (-4,0,0)}
\psfrag{(0,-3,-4)}{\tiny (0,-3,-4)}
\psfrag{(0,-2,-4)}{\tiny (0,-2,-4)}
\psfrag{(0,-1,-5)}{\tiny (0,-1,-5)}
\includegraphics[width=70mm]{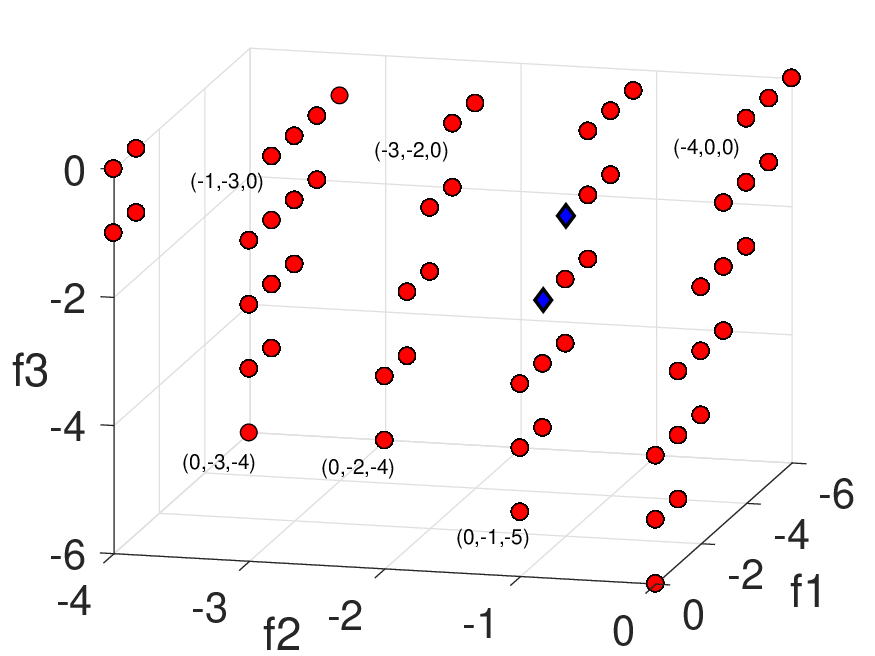} \\
{\footnotesize (c) A rotated view of the front in (a).}
\end{center}
\end{minipage}
\hspace{-1cm}
\begin{minipage}{90mm}
\begin{center}
\hspace*{0cm}
\psfrag{f1}{$f_1$}
\psfrag{f2}{$f_2$}
\psfrag{f3}{\hspace*{-1mm}$f_3$}
\psfrag{(-1,-3,0)}{\tiny (-1,-3,0)}
\psfrag{(-3,-2,0)}{\tiny (-3,-2,0)}
\psfrag{(-4,0,0)}{\tiny (-4,0,0)}
\psfrag{(0,-3,-4)}{\tiny (0,-3,-4)}
\psfrag{(0,-2,-4)}{\tiny(0,-2,-4)}
\psfrag{(0,-1,-5)}{\tiny (0,-1,-5)}
\psfrag{(0,0,-6)}{\tiny (0,0,-6)}
\includegraphics[width=70mm]{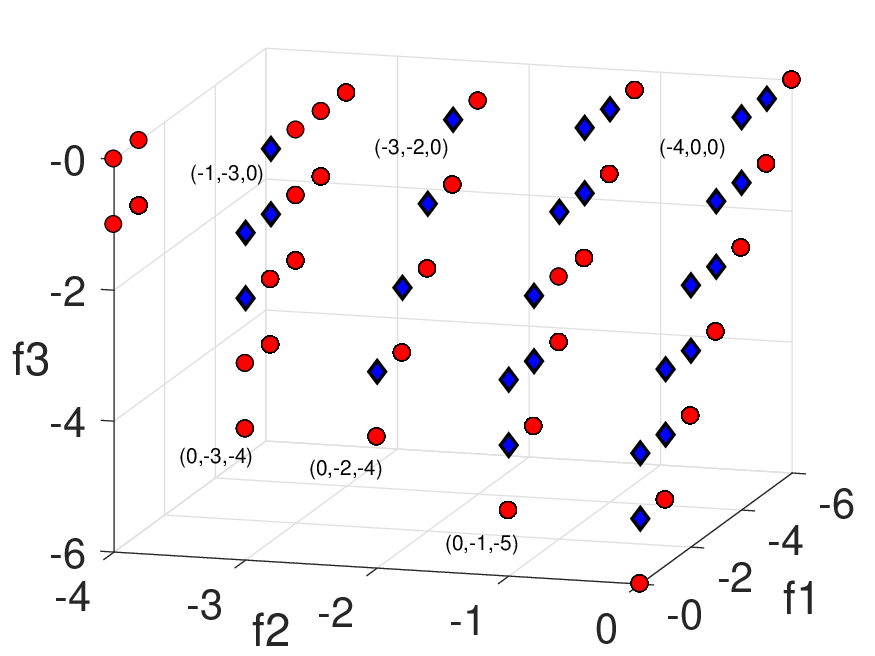} \\
{\footnotesize (d) A rotated view of the front in (b).}
\end{center}
\end{minipage}
\caption{\sf\small Test problem 3 -- Pareto fronts constructed by Algorithms~3 and 4, shown by (red) circles. The (blue) diamonds represent the weak Pareto points that were missed by either algorithm.}
\label{integerAntunes}
\end{figure}

\begin{figure}[hbt!]
\hspace{-1cm}
\begin{minipage}{90mm}
\begin{center}
\hspace*{0cm}
\psfrag{f1}{$f_1$}
\psfrag{f2}{$f_2$}
\psfrag{f3}{\hspace*{-1mm}$f_3$}
\includegraphics[width=75mm]{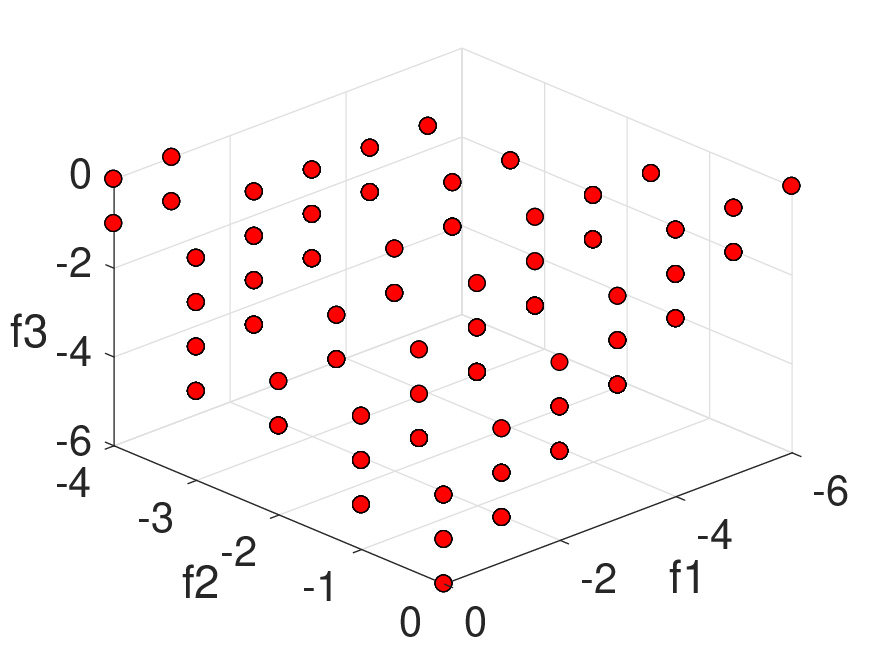} \\
{\footnotesize (a) Pareto front by Algorithm 5.}
\end{center}
\end{minipage}
\hspace{-1cm}
\begin{minipage}{90mm}
\begin{center}
\hspace*{0cm}
\psfrag{f1}{$f_1$}
\psfrag{f2}{$f_2$}
\psfrag{f3}{\hspace*{-1mm}$f_3$}
\includegraphics[width=75mm]{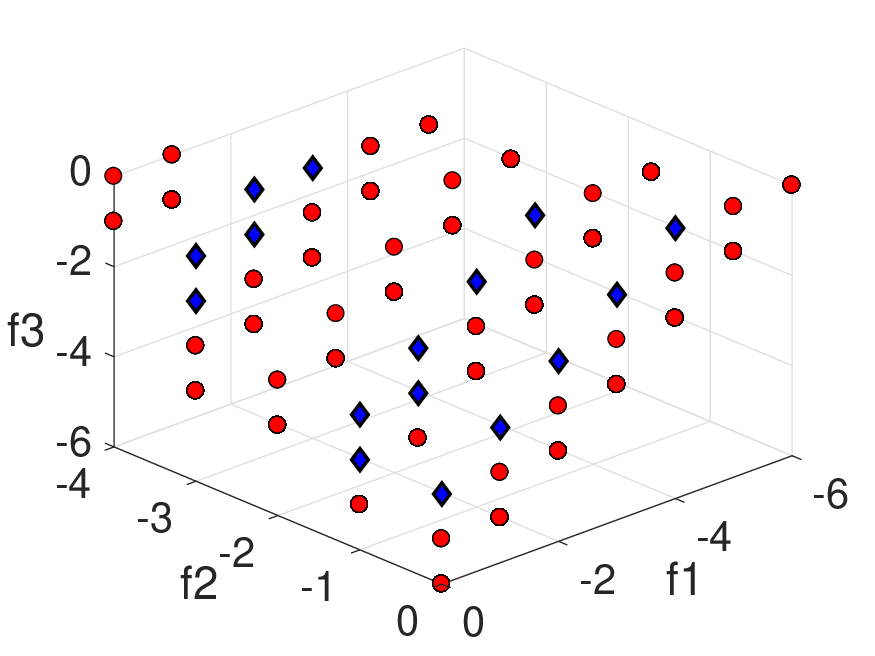} \\
{\footnotesize (b) Pareto front by Algorithm 6.}
\end{center}
\end{minipage}
\\[2mm]
\hspace*{-1cm}
\begin{minipage}{90mm}
\begin{center}
\hspace*{0cm}
\psfrag{f1}{$f_1$}
\psfrag{f2}{$f_2$}
\psfrag{f3}{\hspace*{-1mm}$f_3$}
\psfrag{(-4,0,0)}{\tiny (-4,0,0)}
\psfrag{(0,-3,-4)}{\tiny (0,-3,-4)}
\psfrag{(0,-2,-4)}{\tiny(0,-2,-4)}
\psfrag{(0,-1,-5)}{\tiny (0,-1,-5)}
\includegraphics[width=70mm]{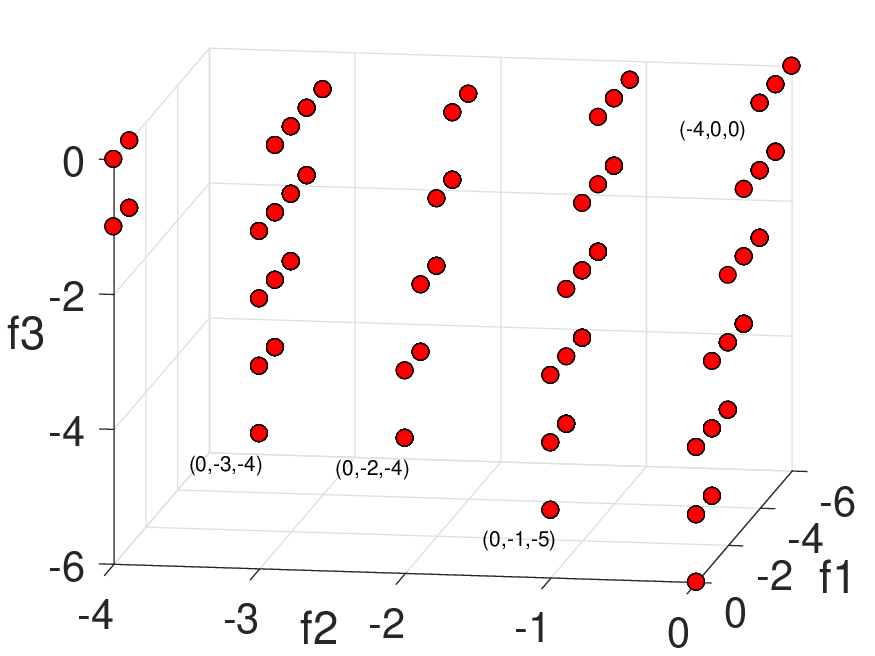} \\
{\footnotesize (c) A rotated view of the front in (a).}
\end{center}
\end{minipage}
\hspace{-1cm}
\begin{minipage}{90mm}
\begin{center}
\hspace*{0cm}
\psfrag{f1}{$f_1$}
\psfrag{f2}{$f_2$}
\psfrag{f3}{\hspace*{-1mm}$f_3$}
\psfrag{(-4,0,0)}{\tiny (-4,0,0)}
\psfrag{(0,-3,-4)}{\tiny (0,-3,-4)}
\psfrag{(0,-2,-4)}{\tiny(0,-2,-4)}
\psfrag{(0,-1,-5)}{\tiny (0,-1,-5)}
\includegraphics[width=70mm]{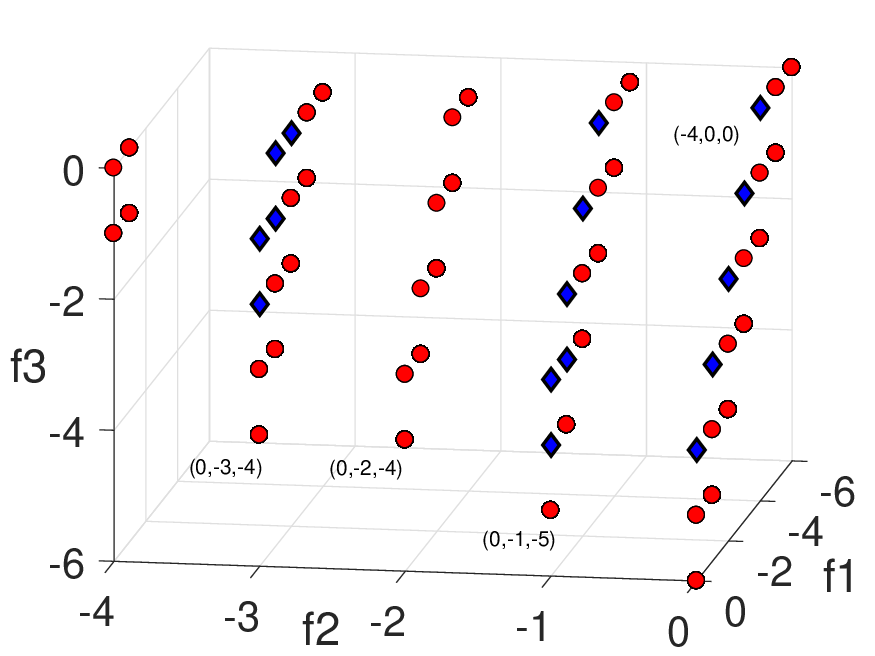} \\
{\footnotesize (d) A rotated view of the front in (b).}
\end{center}
\end{minipage}
\caption{\sf\small Test problem 3 -- Pareto fronts constructed by Algorithms~5 and 6, shown by (red) circles. The (blue) diamonds represent the weak Pareto points that were missed by Algorithm~6.}
\label{testpro1a}
\end{figure}

The three objective functions of the problem are namely $f_1(x)=-x_1$, $f_2(x)=-x_2$ and $f_3(x)=-x_3$. We choose $u = (-100,-100,-100)$ as the utopia point.  For solving this problem, we used Algorithms~3--6. For the features of these algorithms, see Table~\ref{algorithms}.  As in the test problems 1 and 2, for fairness in comparisons, we adjusted the number of points in, depending on the algorithm, the CHIM or SBG grids, in such a way that the CPU time each algorithm takes is roughly the same (around 210 seconds for Algorithms~3--4 and 130 seconds for Algorithms~5--6).  Table~\ref{table:exp-3a} summarizes the numerical performance of Algorithms~3 and 4, while Figures~\ref{integerAntunes}  and \ref{testpro1a} show the Pareto points generated by each algorithm.

For the given CHIM grid, or the allowed CPU time of 202 seconds, Algorithm 3 generates 58 of the 60 weak Pareto points in the front. On the other hand, Algorithm 4 is able to find only around half of the weak Pareto points, in about the same CPU time as that of Algorithm~3.  The Pareto fronts generated by these algorithms are depicted in Figure~\ref{integerAntunes}. When we provide a much finer CHIM grid, i.e., increase the allowed CPU time by about 10 fold, either algorithm fails to generate any new points in the front.  We conclude that Algorithm 3 is more efficient and powerful than Algorithm~4 in approximating the Pareto front.

Algorithm 5, which implements the SBG grid, generates all of the 60 weak Pareto points in the front in just 128 seconds, whereas Algorithm 6  produces only 33 out of the 60 Pareto points, in about the same amount of time. When we provide a much finer SBG grid, resulting in a CPU time of 120 minutes or longer, Algorithm~6 can find at most 44 of the weak Pareto points.  Considering all of Algorithms~3--6, we conclude that Algorithm~5 is the most efficient and powerful.

\begin{table}[t]
\caption{\small{\textit{Test problem 3 -- Numerical performance of Algorithms~3--6.}}}
\footnotesize
\vskip 1.5em
\centering
\begin{tabular}{|c| c| c| c| c|}
\hline
		&  			& Number of 	& Number of  \\
                 & CPU time	& subproblems & weak Pareto  \\
 Algorithm &     [sec] 		&  attempted   	& points generated   \\ [0.5ex]
\hline
3  		& 202 		& $230\times 3 = 690$ & 58\ \ \ \ \ \ \  \\ [.5ex]\hline
4 		& 211 		& 495 		& 33\ \ \ \ \ \ \  \\ [.5ex]\hline
5		& $38+8+82=128$ & $45+91+273=409$  & 60 (all) \\ [.5ex]\hline
6		& $32+14+84=130$ & $60+170+170=400$ & 33\ \ \ \ \ \ \  \\ [.5ex]\hline
6		& $184+1400+5530=7114$ & $360+3081+3081=6522$ & 44\ \ \ \ \ \ \  \\ [.5ex]
\hline
\end{tabular}
\label{table:exp-3a}
\end{table}

To test performance when some of the variables are continuous, we transform Test Problem~3 into a mixed-integer problem in two separate instances:
\begin{enumerate}
    \item [(i)] $x_1$ and $x_2$ are continuous, and $x_3$ an integer.
    \item [(ii)] $x_2$ is continuous, and $x_1$ and $x_3$ integers.
\end{enumerate}

We choose $u = (-100,-100,-100)$ as the utopia point.  To solve either of the instances~(i) and (ii), we use both of Algorithms~3 and 4 and compare.  In order to make comparisons on an equal footing, we adjust the number of CHIM points in the algorithms in such a way that the CPU time each algorithm takes is roughly the same. For case (i), we set around 500 seconds, and for case (ii), around 300 seconds.
\par 

\begin{figure}[hbt!]
\hspace{-1cm}
\begin{minipage}{90mm}
\begin{center}
\hspace*{0cm}
\includegraphics[width=75mm]{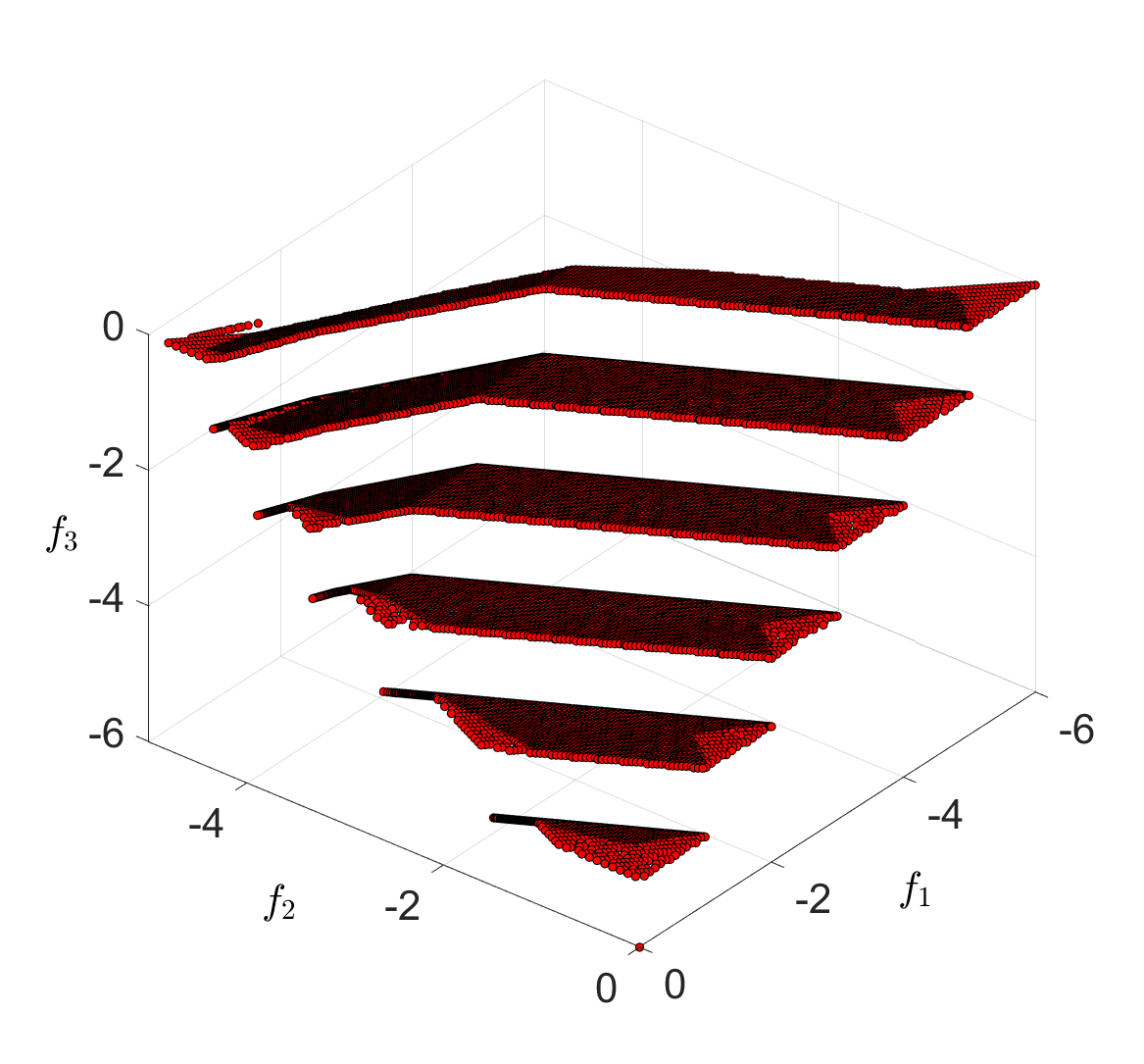} \\
{\footnotesize (a) Instance~(i): Pareto front by Algorithm 3.}
\end{center}
\end{minipage}
\hspace{-1cm}
\begin{minipage}{90mm}
\begin{center}
\hspace*{0cm}
\psfrag{f1}{$f_1$}
\psfrag{f2}{$f_2$}
\psfrag{f3}{\hspace*{-1mm}$f_3$}
\includegraphics[width=75mm]{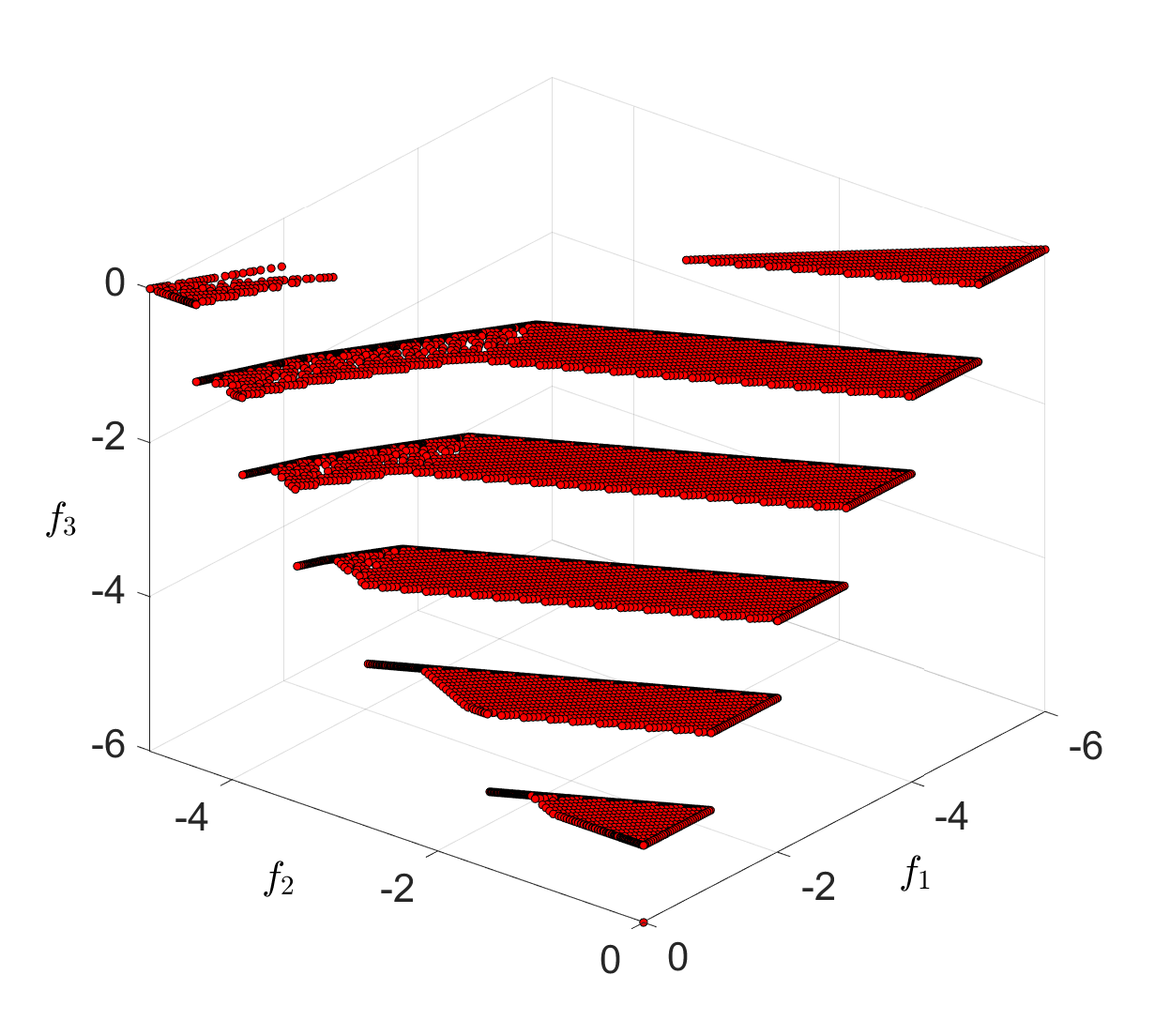} \\
{\footnotesize (b) Instance~(i): Pareto front by Algorithm 4.}
\end{center}
\end{minipage}
\\[2mm]
\hspace*{-1cm}
\begin{minipage}{90mm}
\begin{center}
\hspace*{0cm}
\psfrag{f1}{$f_1$}
\psfrag{f2}{$f_2$}
\psfrag{f3}{\hspace*{-1mm}$f_3$}
\psfrag{(-4,0,0)}{\tiny (-4,0,0)}
\psfrag{(0,-3,-4)}{\tiny (0,-3,-4)}
\psfrag{(0,-2,-4)}{\tiny(0,-2,-4)}
\psfrag{(0,-1,-5)}{\tiny (0,-1,-5)}
\includegraphics[width=70mm]{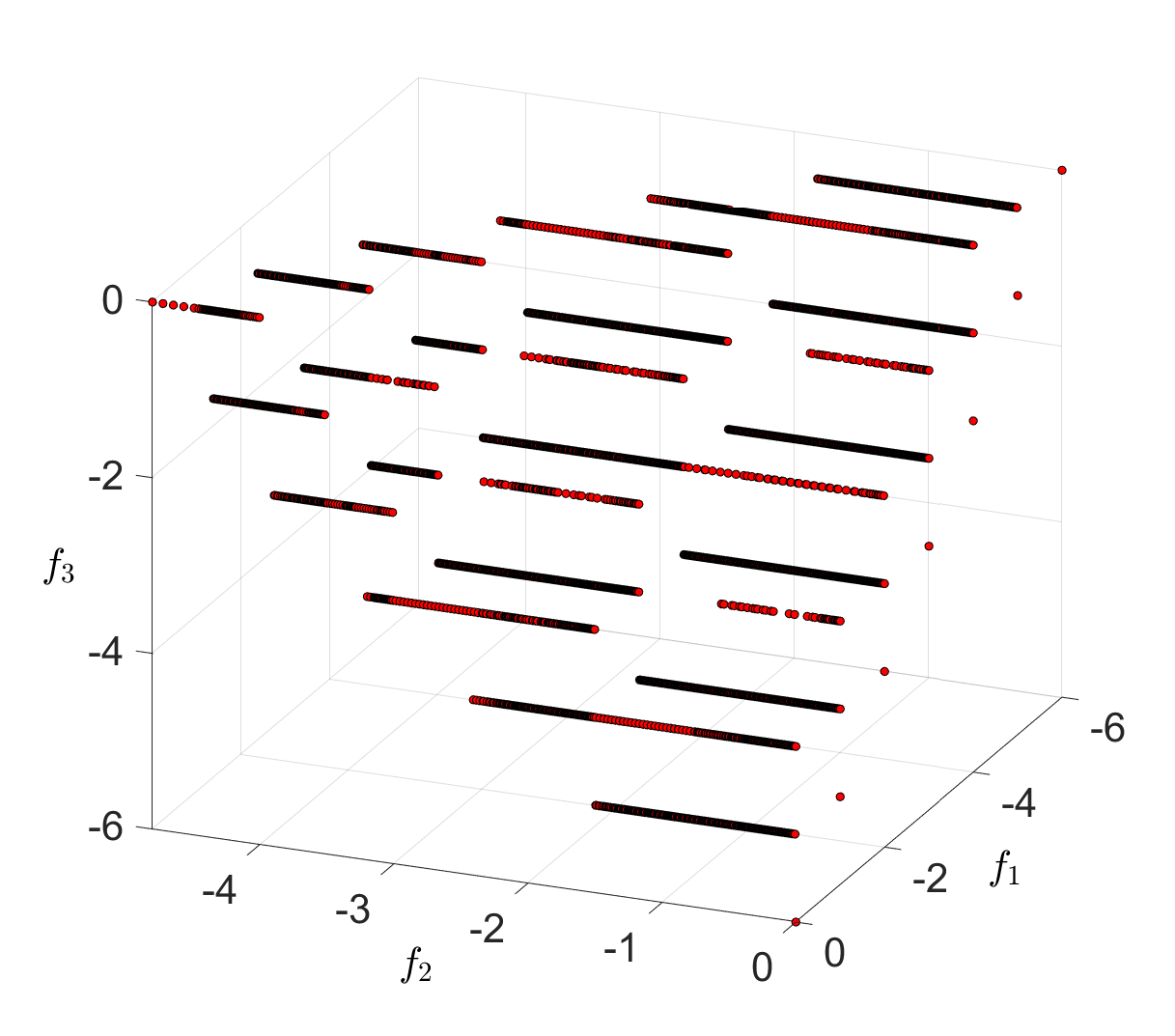} \\
{\footnotesize (c) Instance~(ii): Pareto front by Algorithm 3.}
\end{center}
\end{minipage}
\hspace{-1cm}
\begin{minipage}{90mm}
\begin{center}
\hspace*{0cm}
\psfrag{f1}{$f_1$}
\psfrag{f2}{$f_2$}
\psfrag{f3}{\hspace*{-1mm}$f_3$}
\psfrag{(-4,0,0)}{\tiny (-4,0,0)}
\psfrag{(0,-3,-4)}{\tiny (0,-3,-4)}
\psfrag{(0,-2,-4)}{\tiny(0,-2,-4)}
\psfrag{(0,-1,-5)}{\tiny (0,-1,-5)}
\includegraphics[width=70mm]{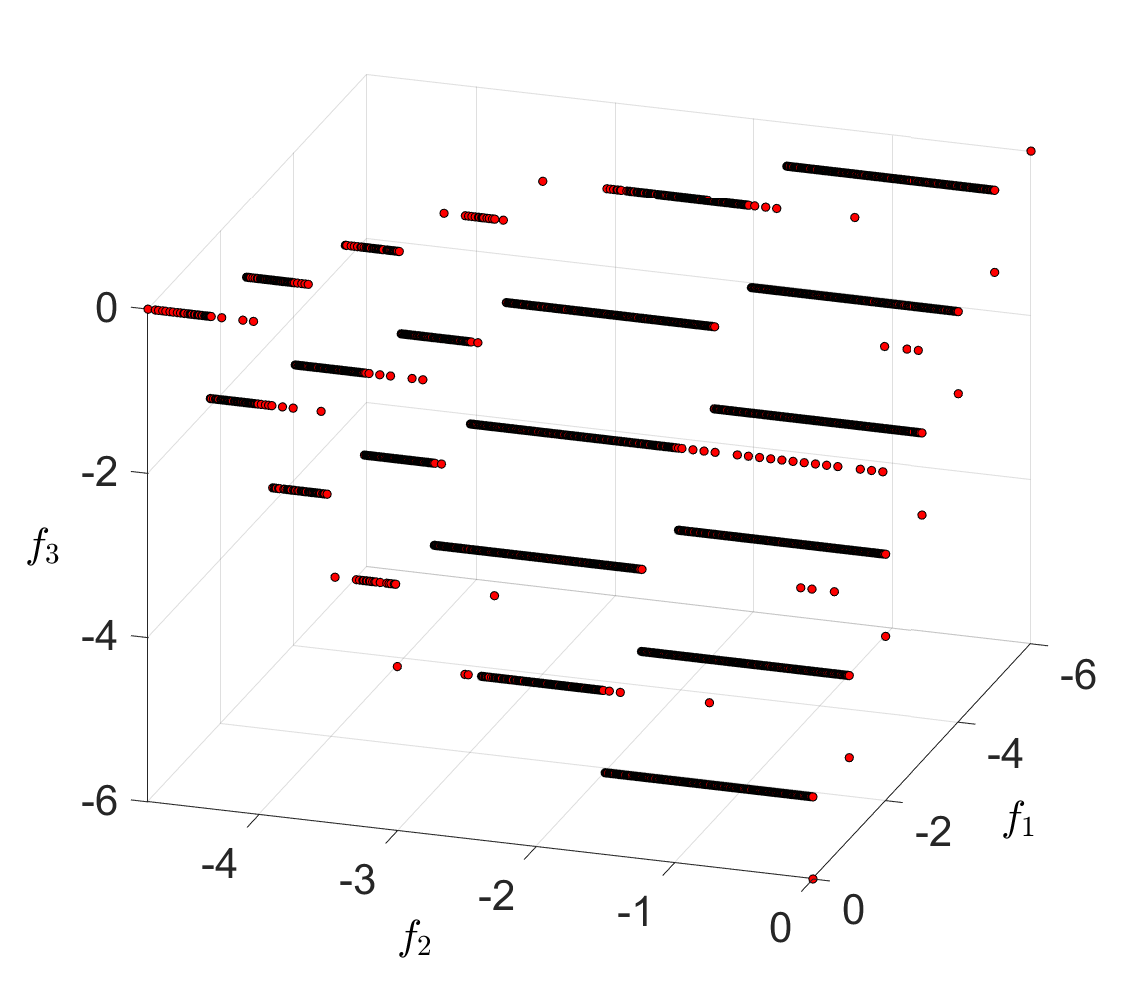} \\
{\footnotesize (d) Instance~(ii): Pareto front by Algorithm 4.}
\end{center}
\end{minipage}
\caption{\sf\small Test problem 3 -- Pareto fronts constructed by Algorithms~3 and 4 for the mixed-integer instances~(i) and (ii).}
\label{testpromix}
\end{figure}

The Pareto fronts obtained for Instance~(i) are shown in Figures~\ref{testpromix}(a) and (b) for Algorithms 3 and 4, respectively. We observe that while Algorithm 4 cannot approximate the upper part of the Pareto front as can be seen in Figure~\ref{testpromix}(b), Algorithm 3 efficiently approximates all six plane segments of the Pareto front as shown in Figure~\ref{testpromix}(a).

In Instance~(ii), we report a similar kind of performance: As can be observed in Figures~\ref{testpromix}(c) and (d), respectively, Algorithm 3 constructs a better approximation of the Pareto front than that by Algorithm~4.  Algorithm~4, which employs the Pascoletti--Serafini scalarization, clearly misses many parts of the front in either instance.

\subsection{Test problem 4: Three objective functions -- nonlinear constraints}
\label{testproblem4}
\begin{figure}[t]
\hspace{-1cm}
\begin{minipage}{90mm}
\begin{center}
\hspace*{0cm}
\psfrag{f1}{$f_1$}
\psfrag{f2}{$f_2$}
\psfrag{f3}{\hspace*{-1mm}$f_3$}
\includegraphics[width=75mm]{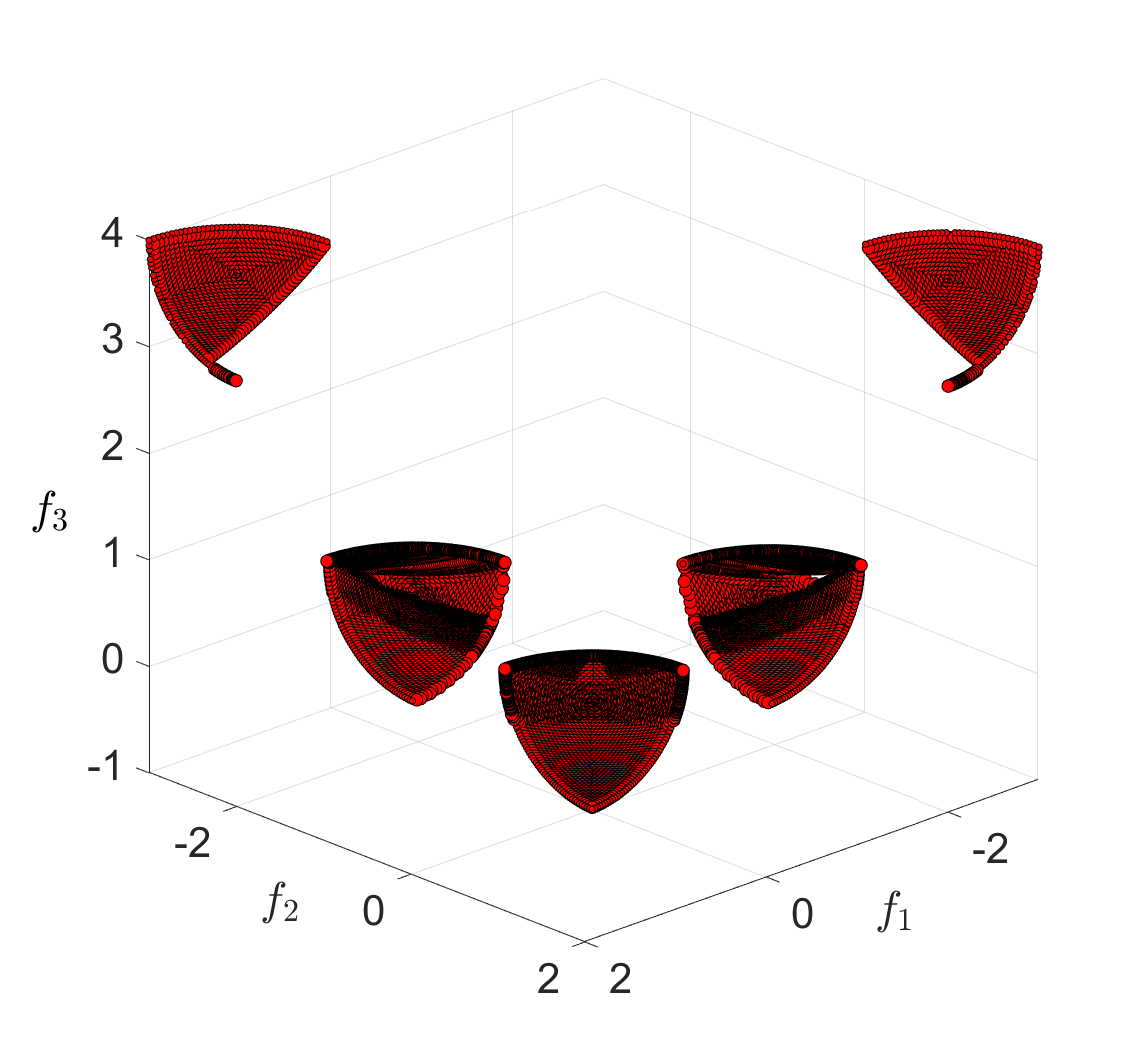} \\
{\footnotesize (a) Pareto front by Algorithm 5.}
\end{center}
\end{minipage}
\hspace{-1cm}
\begin{minipage}{90mm}
\begin{center}
\hspace*{0cm}
\psfrag{f1}{$f_1$}
\psfrag{f2}{$f_2$}
\psfrag{f3}{\hspace*{-1mm}$f_3$}
\includegraphics[width=75mm]{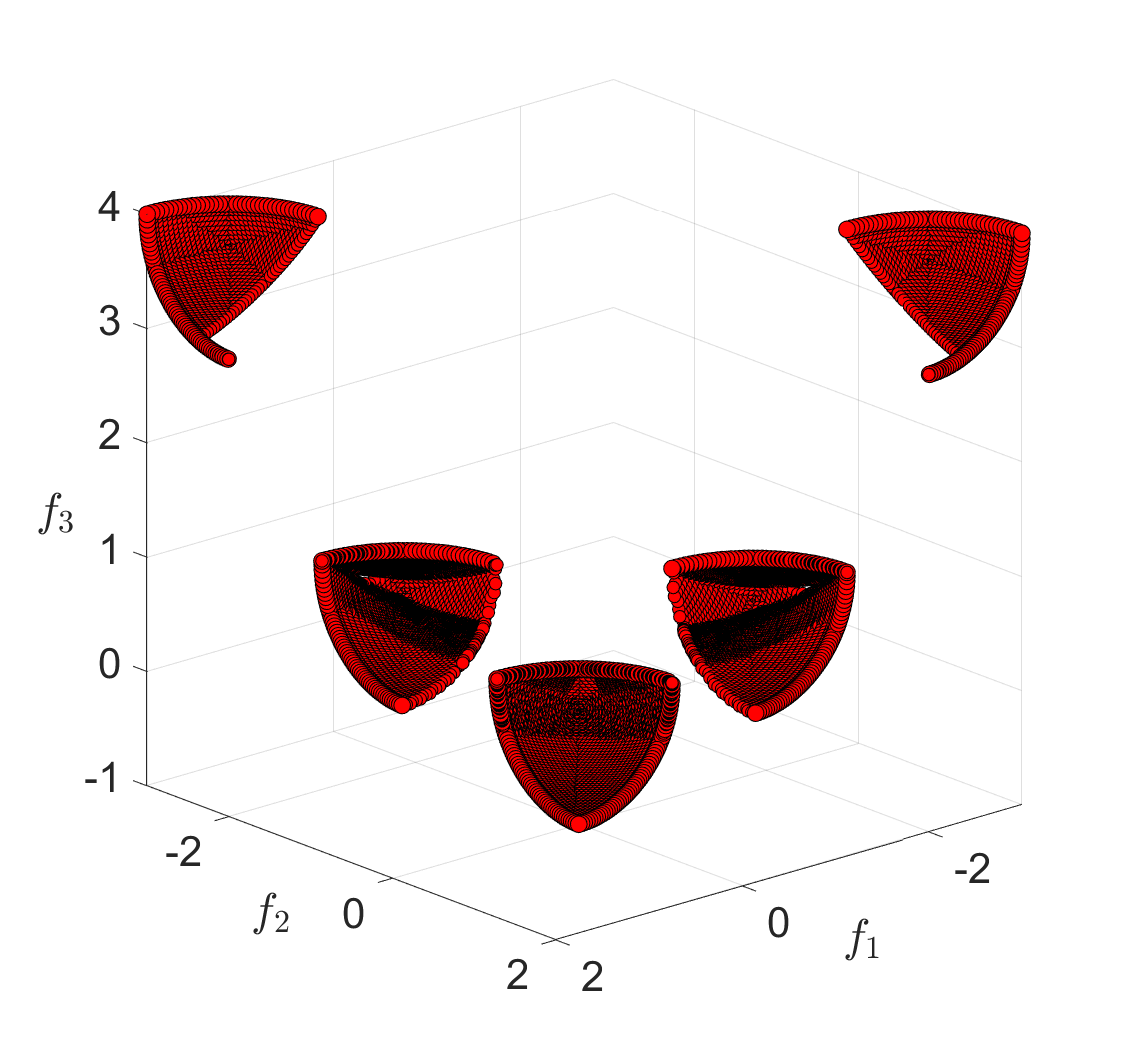} \\
{\footnotesize (b) Pareto front by Algorithm 6.}
\end{center}
\end{minipage}
\hspace*{-1cm}
\begin{minipage}{90mm}
\begin{center}
\hspace*{0cm}
\includegraphics[width=70mm]{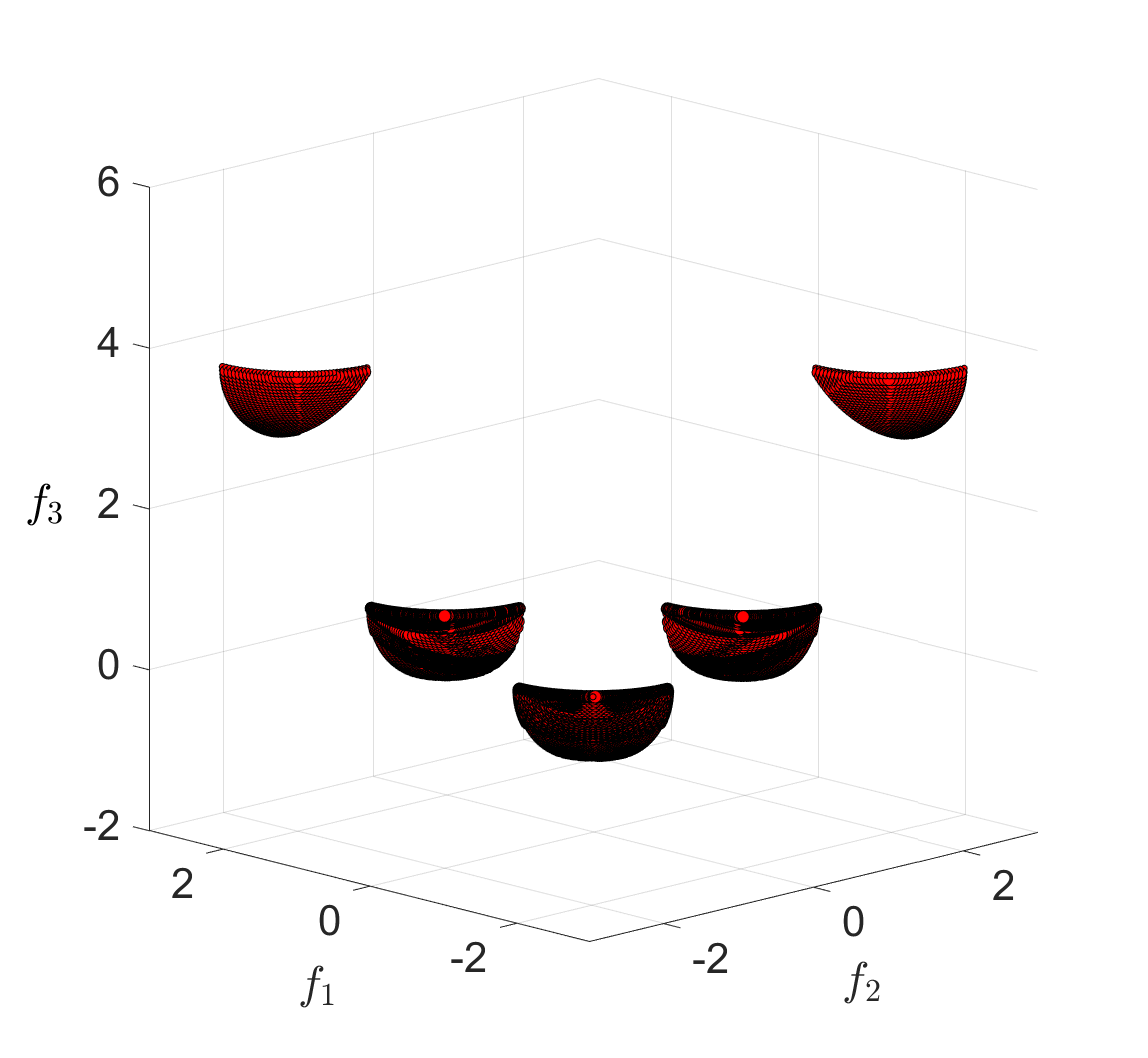} \\
{\footnotesize (c) A rotated view of the front in (a).}
\end{center}
\end{minipage}
\hspace{-1cm}
\begin{minipage}{90mm}
\begin{center}
\hspace*{0cm}
\includegraphics[width=70mm]{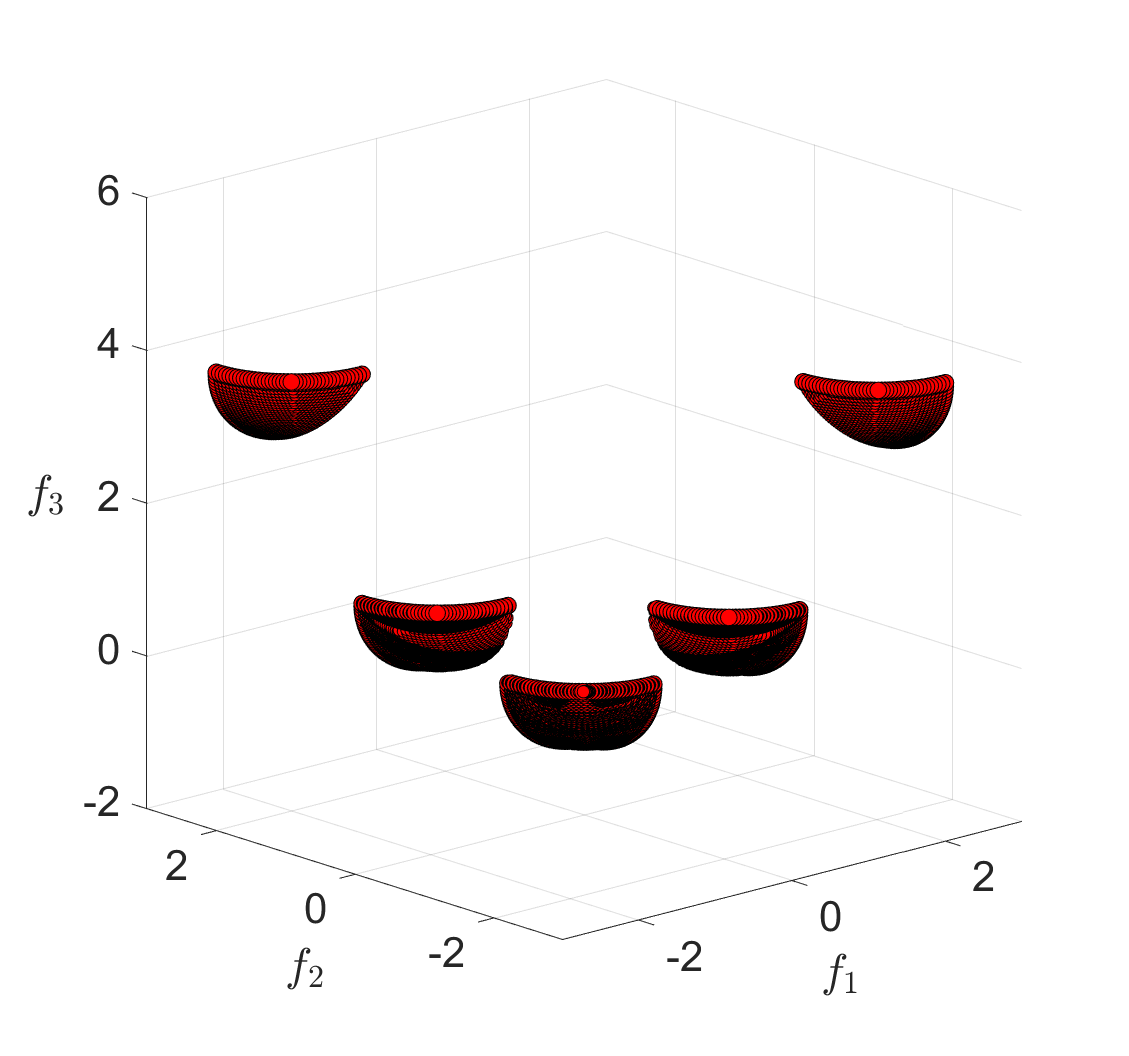} \\
{\footnotesize (d) A rotated view of the front in (b).}
\end{center}
\end{minipage}
\caption{\sf\small Test problem 4 -- Pareto fronts constructed by Algorithms~5 and 6.}
\label{testprogab}
\end{figure}

We consider next a test problem given by de Santis et al.~\cite[Test Instance~4.5, Section~4.3]{Santis2019}, which involves three objective functions and integer, as well as continuous, variables, as stated below.
\[
\begin{array}{cl} \min & \ \left[x_1+x_4,\ x_2-x_4,\ x_3+x_4^2\right] \\[2mm]
 \mbox{s.t} &\  x_1^2+x_2^2+x_3^2 \leq 1, \\[1mm]
 &\ -2\leq x_i\leq 2,\;\; \mbox{for all} \ i=1,\ldots4\,,\ \mbox{ and }  x_4 \mbox{ is an integer}.
\end{array}
\]
The challenge here is to construct (approximately) the whole Pareto front of this problem, and observe the performances of Algorithms 3--6 in carrying out this task.  As can be seen, the problem has three continuous variables and one integer variable.  The resulting Pareto front consists of five disconnected Pareto surfaces---see Figure~\ref{testprogab}. We chose $u = (-100,-100,-100)$ as the utopia point.  We adjusted the number of SBG grid points in the algorithms in such a way that the CPU times needed by the algorithms were more or less the same.  Therefore, their performances could be compared on a more equal footing. For Algorithms 5 and 6, we set around 1200 seconds of CPU time. Both algorithms generated identical Pareto fronts, which are depicted in Figure~\ref{testprogab}.

\begin{figure}[t]
\hspace{-1cm}
\begin{minipage}{90mm}
\begin{center}
\hspace*{0cm}
\psfrag{f1}{$f_1$}
\psfrag{f2}{$f_2$}
\psfrag{f3}{\hspace*{-1mm}$f_3$}
\includegraphics[width=75mm]{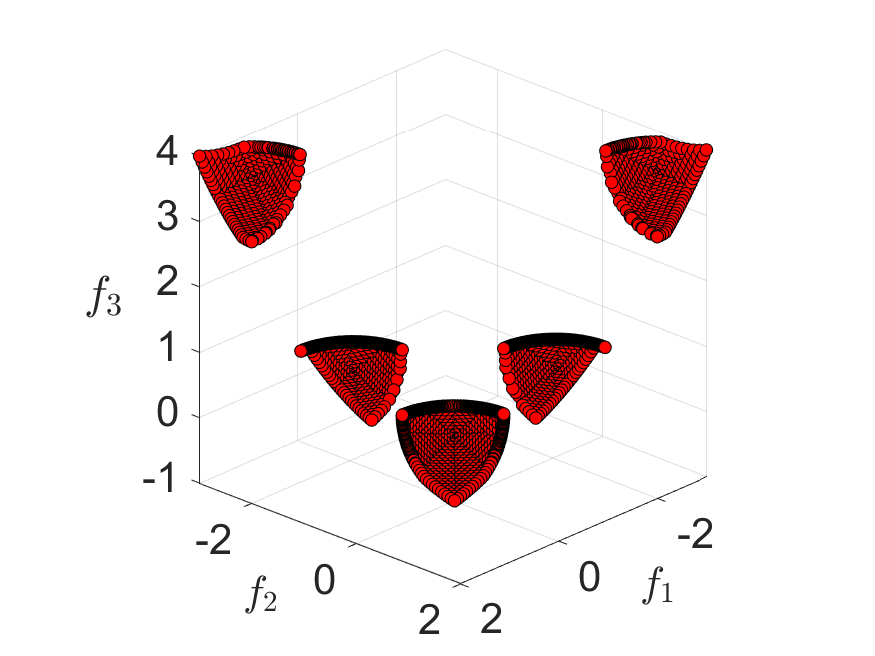} \\
{\footnotesize (a) Pareto front by Algorithm 3.}
\end{center}
\end{minipage}
\hspace{-1cm}
\begin{minipage}{90mm}
\begin{center}
\hspace*{0cm}
\psfrag{f1}{$f_1$}
\psfrag{f2}{$f_2$}
\psfrag{f3}{\hspace*{-1mm}$f_3$}
\includegraphics[width=75mm]{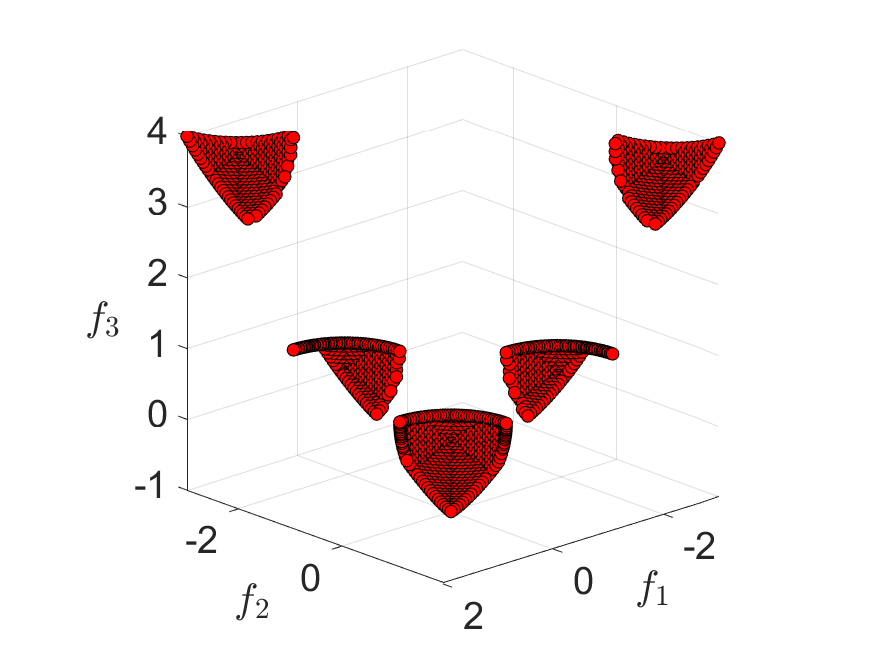} \\
{\footnotesize (b) Pareto front by Algorithm 4.}
\end{center}
\end{minipage}
\caption{\sf\small Test problem 4 -- Pareto fronts constructed by Algorithms~3 and 4.}
\label{testprogabchim}
\end{figure}

We also implemented Algorithms 3 and 4 to solve Test Problem~4.  Although Algorithm~3 can be seen to yield a better approximation of the Pareto front than Algorithm~4  (see Figure~\ref{testprogabchim}), this approximation is still nowhere nearly as good as the one obtained by Algorithms~5--6.  It is conceivable to think that the reason Algorithms~5 and 6 do a better job is that they employ the SBG grid, which is particularly effective when the Pareto front has a complicated boundary, as is the case in this example problem.

\begin{remark}\label{rem:5.1} \rm
 
We aim to generate an approximation of the entire Pareto front, i.e., including also all of the weak efficient points, and we stress that aiming only for the efficient points may miss completely some sizeable and significant portions of the Pareto front. We can elaborate on the significance and justification of this view point as follows.

In the mixed-integer version of Test Problem~3, we find that for any (fixed) assignment of $x_3$ and $x_1$ there is exactly one choice/value of $x_2$ that leads to an efficient solution. By looking at Figure~\ref{testpromix}(c), if we aim only for the efficient points but not the weak efficient ones, any of the line segments appearing in the (weak) Pareto front would be reduced to a single point, leading to a grossly inaccurate approximation of the front.

While our algorithms are successful in obtaining weak efficient solutions which are not efficient, they are also successful in obtaining efficient solutions, when compared with other algorithms. For example, while Algorithm~3 can find all the efficient points in both instances (i) and (ii) of Problem~3---see Figure~\ref{testpromix}(a) and (c)---Algorithm~4 fails in doing so.

\end{remark}

\subsection{An application to rocket injector design}
\label{rocket}

The liquid-rocket injector design problem was previously studied as a multi-objective optimization problem in~\cite{BurKayRiz2017, GoeVaiHafShyQueTuc2007, VaiTucPapShy2004}. Two primary objectives are of concern in this design problem: (i) improvement of the performance of the injector and (ii) increase of its survivability or lifetime. The performance of the injector is influenced by the axial length of the thrust chamber while the lifetime of the injector is associated with the thermal field inside the thrust chamber. For an illustration and visual representation of the injector design we refer the reader to \cite{GoeVaiHafShyQueTuc2007}. There is a conflicting interplay between these two main objectives: while high temperatures improve the performance, they reduce the lifetime and thus decrease survivability. Four design variables are introduced in~\cite{GoeVaiHafShyQueTuc2007} to construct the mathematical model of the rocket injector design problem; namely,
\begin{itemize}
	\item[$x_1$:] hydrogen flow angle,
	\item[$x_2$:] oxidizer post tip thickness,
	\item[$x_3$:] decrement with respect to the baseline cross-section area of the tube carrying oxygen,
	\item[$x_4$:] increment with respect to the baseline cross-section area of the tube carrying hydrogen.
\end{itemize}
In our present study, we introduce an integer-valued variable $\tilde{x}_1\in\{0,1,2,3\}$ and impose a new constraint $x_1=0.2\,\tilde{x}_1$.  The constraint enforces the hydrogen flow angle ($x_1$) to be one of the prescribed set of angles $\{0,0.2,0.4,0.6\}$.  We note that, with the introduction of an integer variable, the modified problem has a disconnected domain and so is expected to have a disconnected Pareto front.

We consider four objective functions as in~\cite{GoeVaiHafShyQueTuc2007}; namely,
\begin{description}
	\item[$f_1(x)$:] face temperature, which is the maximum temperature of the injector face,
	\item[$f_2(x)$:] tip temperature, which is the maximum temperature on the post tip of the injector,
	\item[$f_3(x)$:] combustion length, which is the distance from the inlet where 99\% of the combustion are complete,
	\item[$f_4(x)$:] wall temperature, which is the wall temperature at three inches (fourth probe) from the injector face.
\end{description}
The modified rocket injector design problem considered in this paper is a multi-objective mixed-integer optimization problem and  described as follows.
\[
\min\ \ [f_1, f_2, f_3, f_4]
\]
where
\begin{eqnarray*}
	&& f_1 = 0.692 + 0.477x_1 - 0.687x_4 - 0.08x_3 - 0.065x_2 -
	0.167x_1^2 - 0.0129x_1x_4\\[0mm]
	&&\hspace*{8mm} +\ 0.0796x_4^2 - 0.0634x_1x_3 - 0.0257x_3x_4 +
	0.0877x_3^2 - 0.0521x_1x_2 \\[0mm]
	&&\hspace*{8mm} +\ 0.00156x_2x_4 + 0.00198x_2x_3 + 0.0184x_2^2\,,
	\\[2mm]
	&& f_2 = 0.37 - 0.205x_1 + 0.0307x_4 + 0.108x_3 + 1.019x_2 -
	0.135x_1^2 + 0.0141x_1x_4 \\[0mm]
	&&\hspace*{8mm} +\ 0.0998x_4^2 + 0.208x_1x_3 - 0.0301x_3x_4 -
	0.226x_3^2 + 0.353x_1x_2 - 0.0497x_2x_3 \\[0mm]
	&&\hspace*{8mm} -\ 0.423x_2^2 + 0.202x_1^2x_4 - 0.281x_1^2x_3
	- 0.342x_1x_4^2 - 0.245x_3x_4^2 + 0.281x_3^2x_4 \\[0mm]
	&&\hspace*{8mm} -\ 0.184x_1x_2^2+0.281x_1x_3x_4, \\[2mm]
	&& f_3 = 0.153 - 0.322x_1 + 0.396x_4 + 0.424x_3 + 0.0226x_2 +
	0.175x_1^2 + 0.0185x_1x_4 \\[0mm]
	&&\hspace*{8mm} -\ 0.0701x_4^2 - 0.251x_1x_3 + 0.179x_3x_4 +
	0.015x_3^2  + 0.0134x_1x_2 + 0.0296x_2x_4 \\[0mm]
	&&\hspace*{8mm}  +\ 0.0752x_2x_3 + 0.0192x_2^2, \\[2mm]
	&& f_4 = 0.758 + 0.358x_1 - 0.807x_4 + 0.0925x_3 - 0.0468x_2 -
	0.172x_1^2 + 0.0106x_1x_4 \\ [0mm]
	&&\hspace*{8mm} +\ 0.0697x_4^2 - 0.146x_1x_3 - 0.0416x_3x_4 +
	0.102x_3^2 - 0.0694x_1x_2 \\[0mm]
	&&\hspace*{8mm} -\ 0.00503x_2x_4 + 0.0151x_2x_3+ 0.0173x_2^2\,,
\end{eqnarray*}
subject to
\[
x_1=0.2\,\tilde{x}_1,\;\;
0 \leq \tilde{x}_1 \leq 3\,, \;\; \tilde{x}_1\;\; \mbox{an integer},\]
\[0 \leq x_1, x_2, x_3, x_4 \leq 1\,.
\]

\begin{figure}[t!]
	\hspace{-10mm}
	\begin{minipage}{90mm}
		\begin{center}
			\psfrag{f1}{\small$f_1$}
			\psfrag{f2}{\small$f_2$}
			\psfrag{f3}{\small$f_3$}
			\includegraphics[width=85mm]{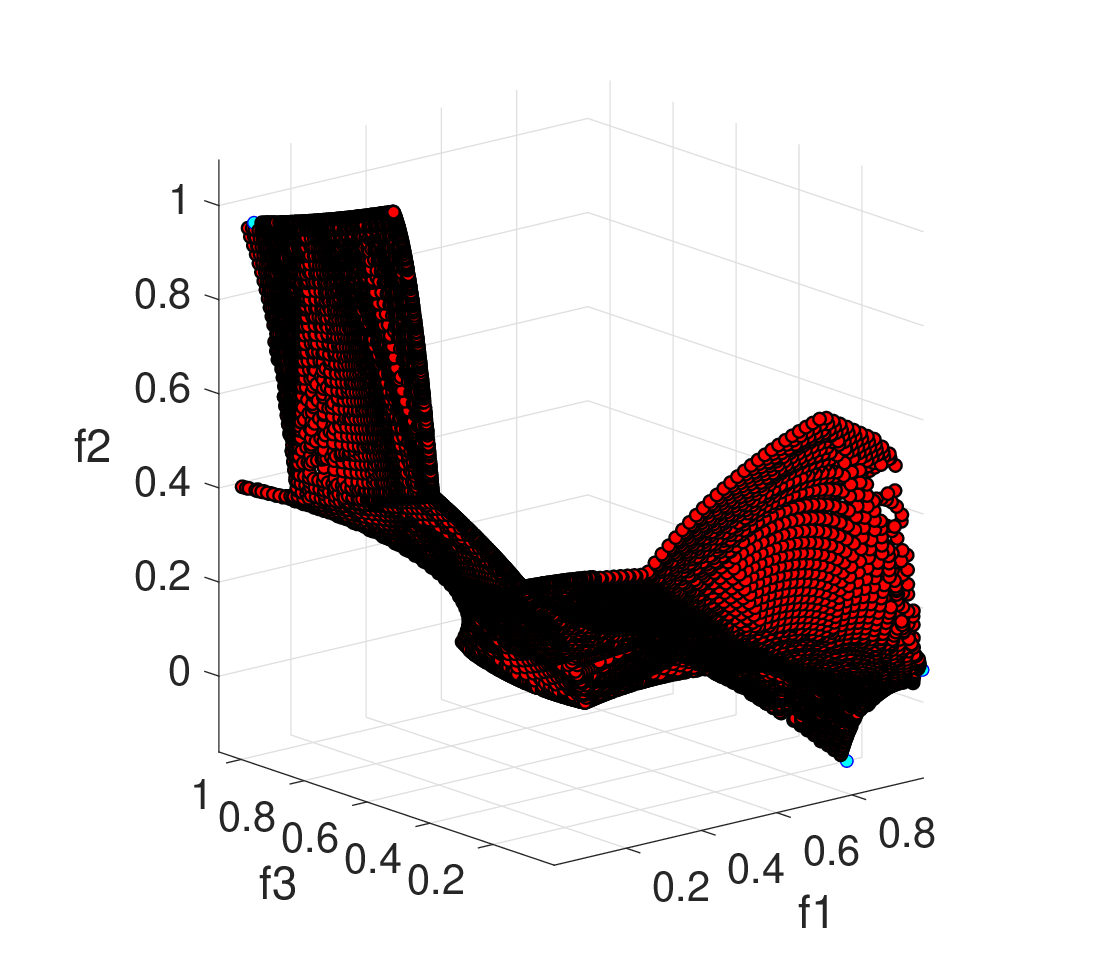} \\
			{\scriptsize (a) Projection of the Pareto front\\ to the $f_1f_3f_2$-space \cite{BurKayRiz2017}.}
		\end{center}
	\end{minipage}
	\hspace{0mm}
	\begin{minipage}{90mm}
		\begin{center}
		\hspace{-20mm}
			\psfrag{f1}{\small$f_1$}
			\psfrag{f2}{\small$f_2$}
			\psfrag{f3}{\small$f_3$}
			\includegraphics[width=85mm]{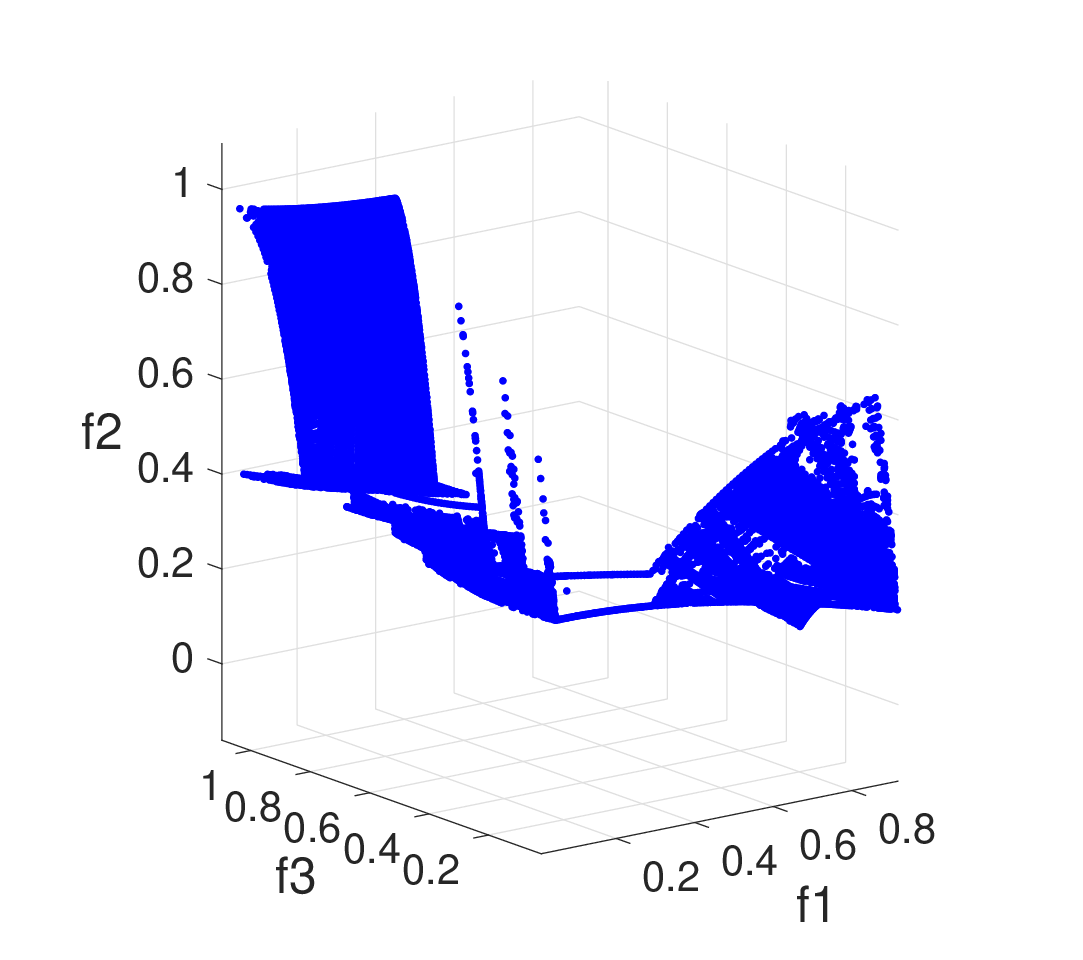} \\
			{\scriptsize \hspace*{-20mm} (b) Projection of the Pareto front\\
				\hspace*{-20mm} generated by Algorithm~7.}
		\end{center}
	\end{minipage}
\\[2mm]
\hspace*{-1cm}
\begin{minipage}{90mm}
	\begin{center}
		\hspace*{0cm}
		\psfrag{f1}{$f_1$}
		\psfrag{f2}{$f_2$}
		\psfrag{f3}{$f_3$}
		\includegraphics[width=85mm]{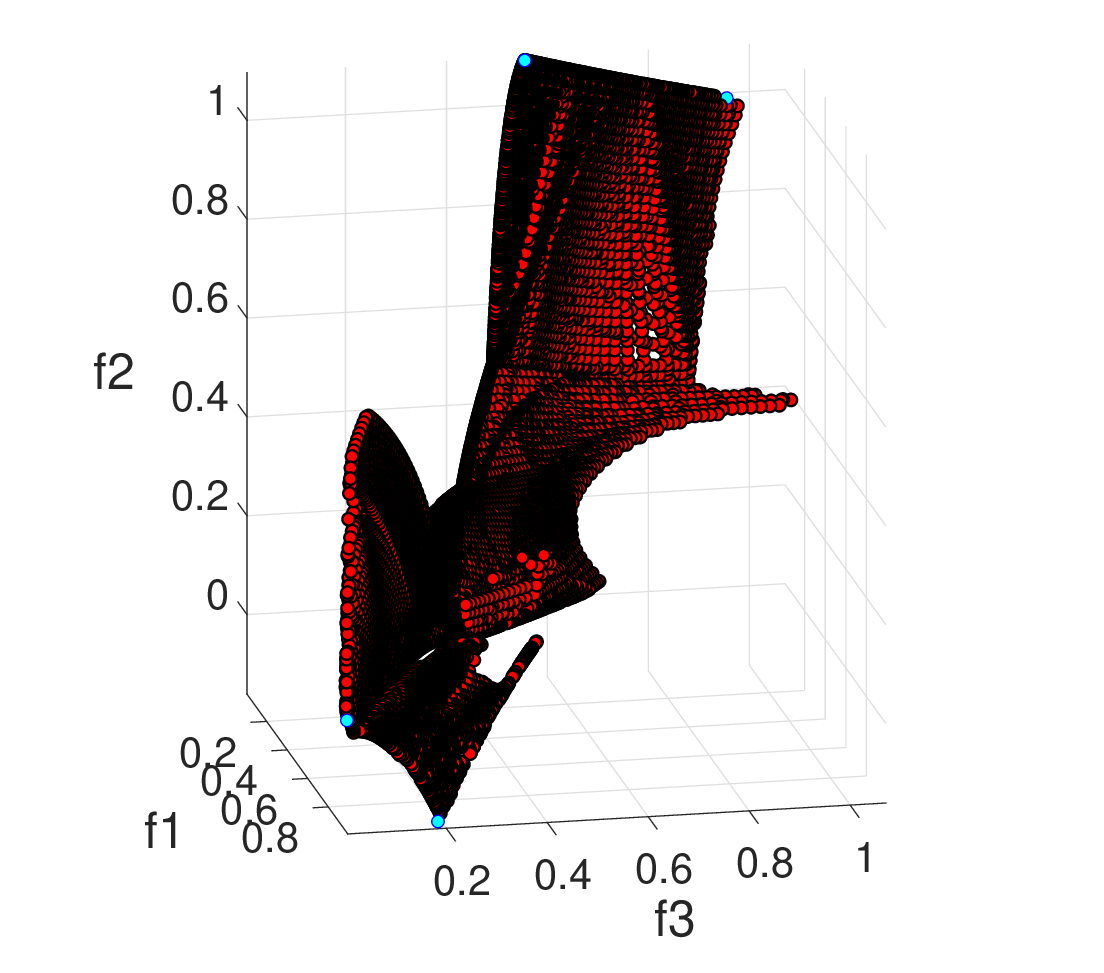} \\
		{\scriptsize (c) Rotated view of the front in (a) \cite{BurKayRiz2017}.}
	\end{center}
\end{minipage}
\hspace{-0cm}
\begin{minipage}{90mm}
	\begin{center}
		\hspace*{-2cm}
		\psfrag{f1}{$f_1$}
		\psfrag{f2}{$f_2$}
		\psfrag{f3}{$f_3$}
		\includegraphics[width=100mm]{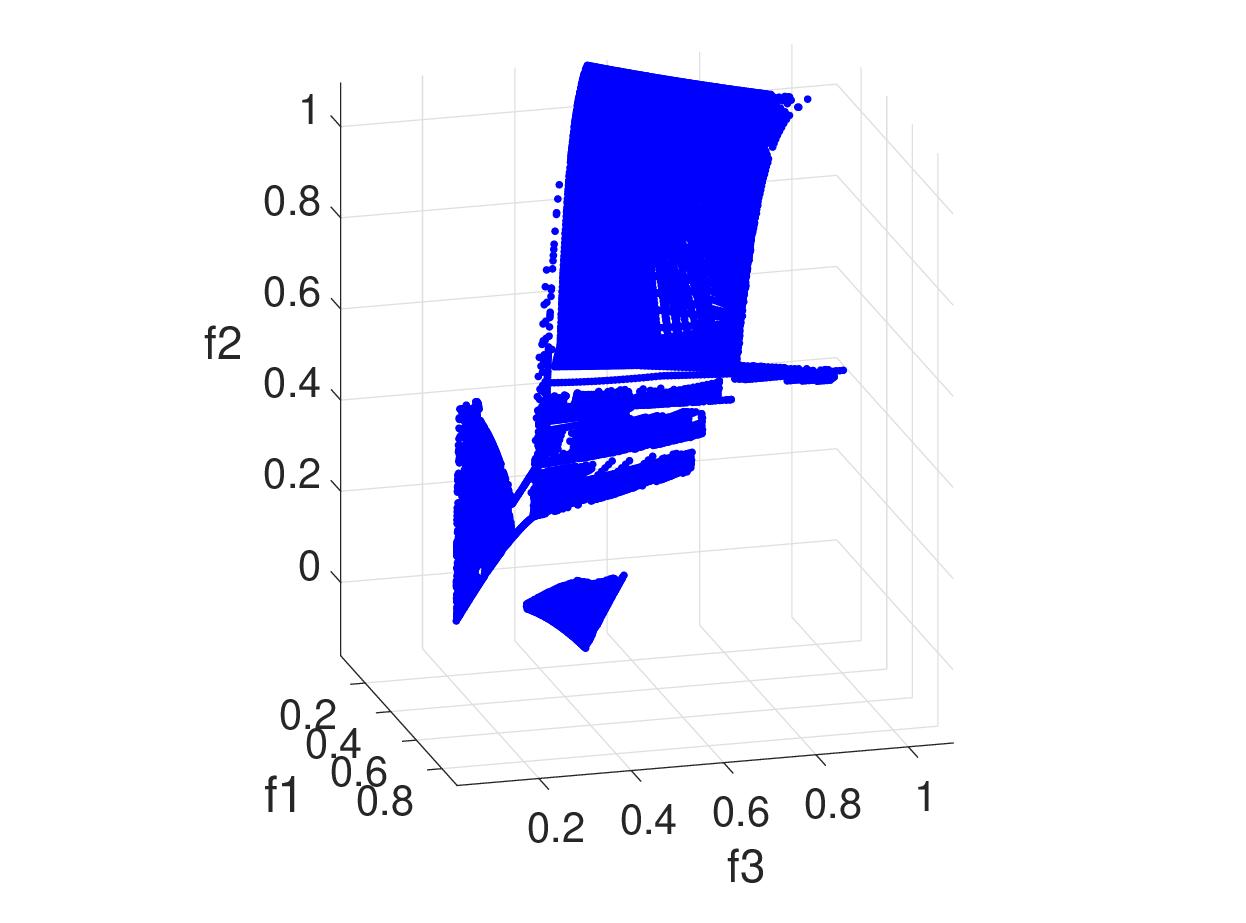} \\
		{\scriptsize (d) Rotated view of the front in (b).}
	\end{center}
\end{minipage}
	\caption{Projected Pareto fronts for the rocket injector design problem. For comparison purposes, the case of continuous variables \cite{BurKayRiz2017} is depicted in Figures~5(a) and (c).}
	\label{injector}
\end{figure}
In obtaining an approximation of the Pareto front of the above problem in the four-dimensional value space, we have used  the new Algorithm~7, which implements the weighted-constraint scalarization and the SBG grid.  Note that an implementation of the CHIM grid generation technique (instead of the SBG grid) in Algorithm~7 would not have worked for the rocket injector design problem, because of the rather complex boundary of the Pareto front of this problem. Hence Algorithm~7 uses the SBG grid.  

We compare our mixed-integer multi-objective programming results with those obtained for the case of continuous variables in \cite{BurKayRiz2017}. We choose to display the projection $(f_1,f_2,f_3,f_4)\mapsto (f_1,f_2,f_3,0)$ in the $f_1f_3f_2$-space for this comparison.  The Pareto front depicted in parts~(a) and (c) of Figure~\ref{injector} corresponds to the case of continuous variables in \cite{BurKayRiz2017}, and was obtained by using Algorithm~9 from \cite{BurKayRiz2017}.  On the other hand, the Pareto front in parts~(b) and (d) of Figure~\ref{injector} corresponds to the case of mixed integer variables, and have been obtained by using Algorithm~7.

Algorithm~7 obtained $41014$ points, which, after the elimination of dominated points, resulted in $ 33536$ points representing the Pareto front. Note that what we obtain is actually a discretization of the front, i.e., the continuous parts of the front are approximated by, in particular computed as, a discrete set of points.  A computed point may in fact not be a weak Pareto point due to computational errors.  Therefore, it is possible to classify such erroneous points as non-dominated when compared with the other computed points.  For this reason, one cannot ensure that all dominated points have been eliminated, or weeded out; however, this example application still illustrates the utility of our approach.  The Pareto front approximation in Figure~\ref{injector}(a) for continuous variables had taken about 2 hours, using the same computer, as reported in \cite{BurKayRiz2017}.  The Pareto approximation in Figure~\ref{injector}(b), on the other hand, took about 32 hours. One should recall that solving mixed integer problems in general requires  substantially more computer memory and CPU time.

\subsubsection{Comparison of Algorithm 7 with a brute force or slicing technique}
In order to obtain an approximation of the Pareto front, we aim to propose algorithms supported with mathematical results rather than a {\em brute-force} methodology.   The following reasons motivate our approach for not resorting to brute-force.
\begin{itemize}
    \item  The computational time required by a brute-force algorithm will depend exponentially on the number of optimization variables and the number of objective functions.
    \item  A more accurate representation of the front would typically require an exponentially growing number of points and effort, which might grind the process down to a halt.
    \item  If the Pareto front's boundary is fragmented or disconnected, i.e., if the front itself has gaps on the ``inside of itself," then one would need even more points to approximate the Pareto front. Many of the points generated by algorithms around these gaps are typically non-Pareto points and, therefore, require more computational time.  This case creates a bigger concern for the brute-force algorithm.
\end{itemize}

Here we aim to justify our choice of not resorting to brute force, or the so-called {\em slicing} technique, by carrying out experiments on the rocket injector design problem with each of the methods (brute and non-brute) and making comparisons.  Both methods are allowed to take about (the same) 2.8 hours of CPU time.  Here we first run Algorithm 7 for 2.8 hours with slicing, i.e., run Algorithm 7 four times consecutively, each time with a fixed value of $x_1 = 0, 0.2, 0.4, 0.6$, and superpose all the solutions and do the necessary weeding (or filtering) out.  We run Algorithm 7 again for 2.8 hours, but this time only once, without fixing the value of $x_1$.  Ultimately we compare the performances of each approach by looking the Pareto fronts that they generate and see which of these Pareto fronts come closer to the refined/accurate front presented in Figures~\ref{injector}(b) and \ref{injector}(d).

In the approximation of the Pareto front depicted in Figure~\ref{injector}(b) and \ref{injector}(d), we had used the mixed-integer problem solver SCIP for the subproblems of Algorithm 7.  On the other hand, when we apply the slicing technique we get a continuous optimization problem for every slice, i.e., for every fixed value of $x_1$. Therefore we use a continuous solver, namely Ipopt with a MATLAB interface, for the slicing technique, in order to exploit fully the computational power of a continuous solver.  It is well-known that continuous solvers are far more efficient than mixed-integer solvers in solving a continuous problem, in general.

In running Algorithm~7 for the rocket injector design problem, we have adjusted the number of SBG grid points in Algorithm 7 so as to adjust the CPU time to be about $2.8$ hours. For each of the four fixed values of $x_1$, Algorithm 7 has generated four Pareto fronts for the four {\em sliced-up} problems. Adding all the points in these fronts has then yielded $21611$ ``candidate'' solutions. Discarding/weeding out the dominated points leaves over a mere $8571$ non-dominated points, which are shown in Figures~\ref{disrock}(a) and \ref{disrock}(c). This illustrates that only about 40\% of the points found by Algorithm~7 via slicing are non-dominated.  

We have also solved the mixed-integer rocket injector design problem by running Algorithm~7 without any slicing, i.e., without fixing the value of $x_1$. In this case, we have set the same 2.8 hours of CPU time by providing the size of the SBG grid to Algorithm 7 accordingly. Since the subproblems of Algorithm 7 (without slicing) are mixed-integer problems, we have tried each of the mixed-integer solvers Bonmin, BARON and SCIP, and observed that Bonmin was the fastest. With the choice of Bonmin for solving the subproblems, Algorithm 7 generated $21882$ candidate points.  Of these points, after weeding out, we obtained $19698$ non-dominated points, which are about 90\% of all points computed. In this particular example, the number of Pareto points found by Algorithm~7 without slicing is about 2.3 times higher than that with slicing.

The Pareto front obtained by Algorithm 7 without slicing is depicted in Figures~\ref{disrock}(b) and \ref{disrock}(d) from two distinct viewing angles.  These viewing angles are the same as those used in Figures~\ref{injector}(b) and \ref{injector}(d), and in \ref{disrock}(a) and \ref{disrock}(c), respectively, for ease of comparisons. First, we observe that a triangle-shaped part of the front is missing from the approximation of the front presented in Figure~\ref{disrock}(c)---compare with Figures~\ref{injector}(d) and \ref{disrock}(d). We also observe that the boundary of the front in Figures~\ref{disrock}(b) and \ref{disrock}(d) look better developed and defined.  Overall, the images of the front in Figures~\ref{disrock}(b) and \ref{disrock}(d) seem to be closer in appearance to those in Figures~\ref{injector}(b) and \ref{injector}(d), compared with those in Figures~\ref{disrock}(a) and \ref{disrock}(c).  Therefore, based on the numerical experiments and results for this particular example, we conclude that Algorithm~7 without slicing does a better job than that with slicing.

We note that filtering, or weeding out, was also carried out for all the test examples other than the rocket injector design problem, but no dominated points were found in those cases. Filtering seems to be more important, or necessary, for the rocket injector design problem, since this problem is much more difficult than the previous ones. There are cases in the rocket injector design problem when a point computed by the optimization software is only locally Pareto optimum. In such a case, we need to weed out these points, if we can, by making comparisons with the other points computed. However, if one has not found many other Pareto points around a ``local'' (dominated) Pareto point, then it may not be possible to remove that dominated only-locally-Pareto point.  We have to accept, on the other hand, that this situation is an inherent deficiency/property of any general algorithm, unless the problem in question has a particular structure that can be exploited to find global Pareto optimal points.
\textcolor{blue}{
\begin{figure}[t!]
\hspace{-1cm}
\begin{minipage}{90mm}
\begin{center}
\hspace*{-1cm}
\psfrag{f1}{$f_1$}
\psfrag{f2}{$f_2$}
\psfrag{f3}{\hspace*{-1mm}$f_3$}
\includegraphics[width=105mm]{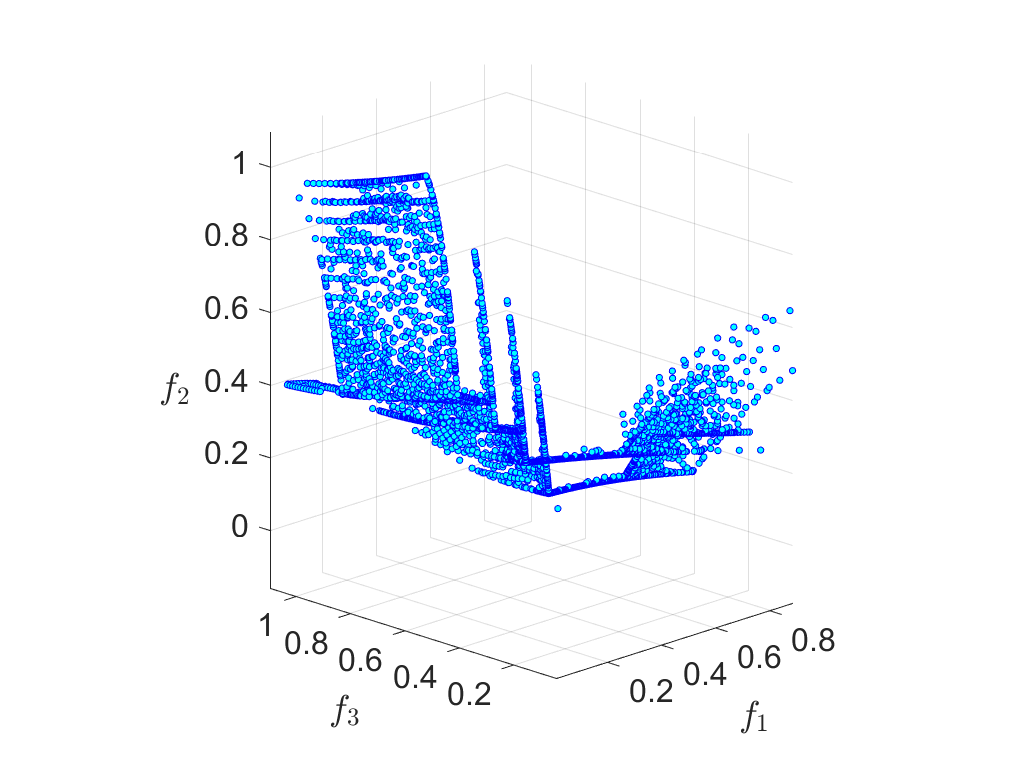} \\
{\footnotesize (a) Projected Pareto front by Algorithm~7 with slicing.}
\end{center}
\end{minipage}
\hspace{-1cm}
\begin{minipage}{90mm}
\begin{center}
\hspace*{0cm}
\psfrag{f1}{$f_1$}
\psfrag{f2}{$f_2$}
\psfrag{f3}{\hspace*{-1mm}$f_3$}
\includegraphics[width=95mm]{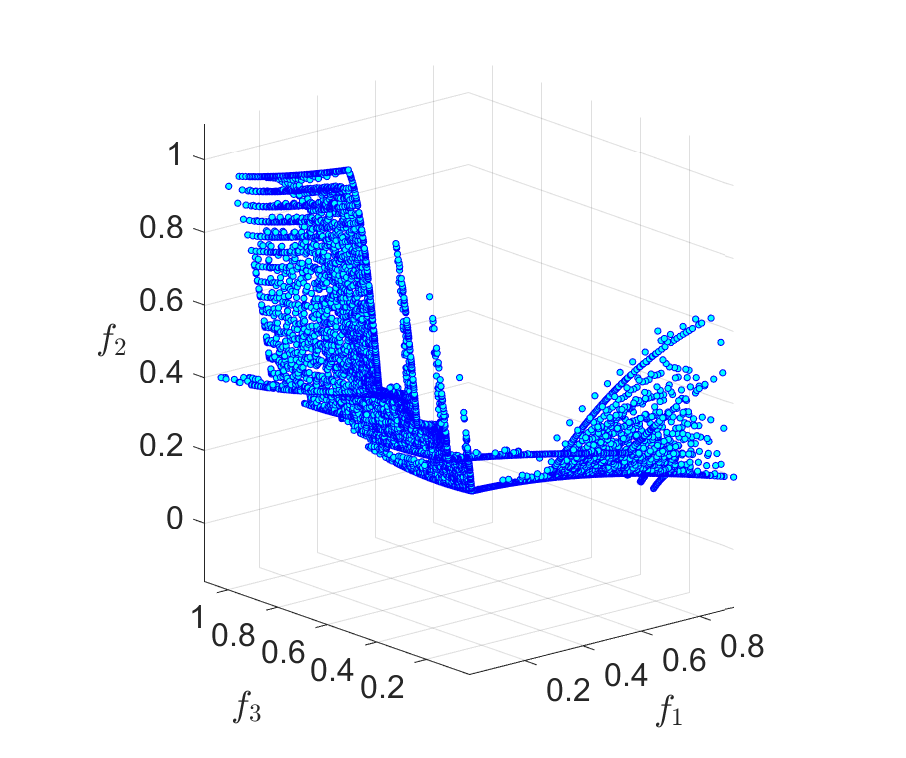} \\
{\footnotesize (b) Projected Pareto front by Algorithm~7 without slicing.}
\end{center}
\end{minipage}
\hspace*{-1.5cm}
\begin{minipage}{90mm}
\begin{center}
\hspace*{0cm}
\includegraphics[width=95mm]{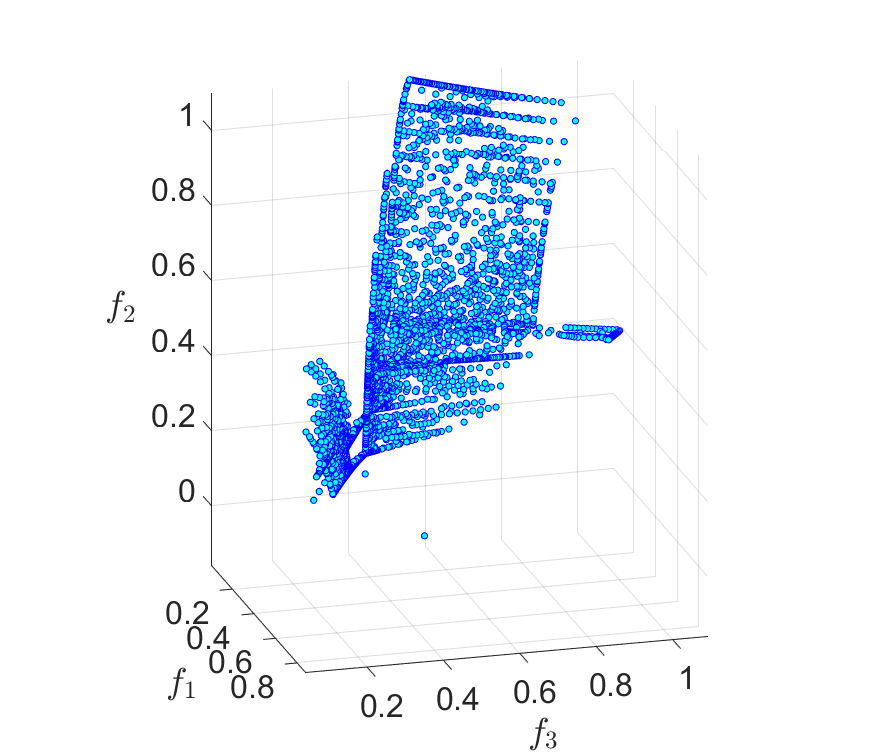} \\
{\footnotesize (c) A rotated view of the front in (a).}
\end{center}
\end{minipage}
\hspace{-1cm}
\begin{minipage}{90mm}
\begin{center}
\hspace*{0cm}
\includegraphics[width=95mm]{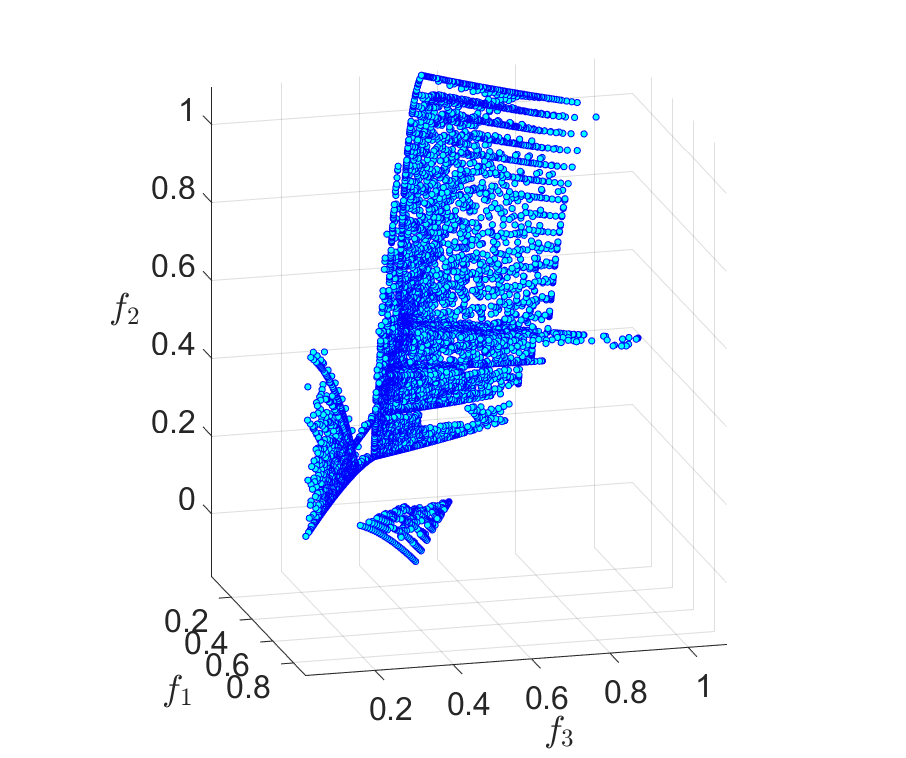} \\
{\footnotesize (d) A rotated view of the front in (b).}
\end{center}
\end{minipage}
\caption{\sf\small  Rocket injector design problem -- Pareto fronts constructed by Algorithm 7, with slicing (parts (a) and (c)) and without slicing (parts (b) and (d)), after 2.8 hours of run.}
\label{disrock}
\end{figure}
}

\section{Conclusion and Discussion}\label{sec6}
\label{conc}

We have established relationships amongst various scalarization techniques listed in Section~\ref{sec3}. We have implemented and tested algorithms, which were originally proposed for continuous multi-objective problems, for constructing the Pareto front of integer and mixed-integer programming problems.  These algorithms utilize one of the CHIM and SBG grid generation techniques and one of the Pascoletti--Serafini or weighted-constraint scalarization techniques, with two or three objective functions to minimize.  We have also presented a new four-objective algorithm.  We applied the algorithms to test problems, one of them being a challenging four objective mixed-integer programming problem.  

We found that for the given test problems, (i)~Algorithm~1 outperforms Algorithm~2 for the two-objective test problem, (ii)~Algorithm~3 outperforms Algorithm~4 for the first three-objective test problem and (iii)~Algorithm~5 outperforms Algorithm~6 for the second three-objective test problem.  The common characteristic of the algorithms which are successful is that they use the weighted-constraint scalarization.  This scalarization approach was previously demonstrated in \cite{BurKayRiz2014, BurKayRiz2017} to perform well for continuous-variable problems with disconnected feasible domains.  It is also worth noting that, in solving the second three-objective test problem, we have demonstrated that Algorithm~5, which utilizes the SBG grid and the weighted-constraint scalarization, outperforms not only Algorithm~6 but also Algorithms~3 and 4.

We have used the weighted-constraint scalarization and the SBG grid in Algorithm~7, which we have presented in this paper for the first time -- see Appendix~A.2.  We have proposed a modified, mixed-integer programming, version of the already challenging rocket injector design problem (previously studied as a continuous problem in~\cite{BurKayRiz2017}), and applied Algorithm~7 to construct the Pareto front of this problem successfully.  The choice of the SBG grid for Algorithm~7 is essential as the boundary of the Pareto front of this problem is quite complicated.  

The continuous version of the rocket injector design problem was solved in~\cite{BurKayRiz2017} with an algorithm (Algorithm~9 in~\cite{BurKayRiz2017}) implementing the feasible-value-constraint scalarization (also introduced in~\cite{BurKayRiz2017}) and the SBG grid.  In the same paper, an algorithm with the same scalarization and the CHIM grid (instead of the SBG grid) was not successful.  In view of Remark~\ref{rem:choice} and the particular success of the weighted-constraint scalarization with problems with disconnected domains, we have proposed in the current paper Algorithm~7 implementing this scalarization along with the SBG grid.

\section*{Acknowledgments}

The authors would like to thank J{\"o}rg Fliege once again for pointing to the challenging multi-objective rocket injector design problem, which was earlier used by the authors in \cite{BurKayRiz2017}. In the present paper, the rocket injector problem has been modified/altered so as to make it a multi-objective mixed-integer problem, yielding an even more challenging instance.

\section*{Appendix A}

\subsection*{A.1\ \ Algorithm 3}

We use the CHIM grid and the weighted-constraint scalarization in the algorithm.  In Step 2 of Algorithm 3, each objective function, subject to the original constraints of the problem, is minimized. These individual minima are used in Step 3 to form a triangle shaped grid. Each grid point corresponds to a weight vector in $\mathbb{R}^3$. In Step 4, three sub-problems are solved at each grid point to generate Pareto points. \par
In Step~4(b), the approximation of the weak efficient point is calculated using the fact that if $\bar{x}_1=\bar{x}_2=\bar{x}_3=:\bar{x}\; \mbox{(say)}$ holds, then the solution $\bar{x}$ would be weak efficient. Here,  $\bar{x}_k$ are the solutions of (P$_w^k$) for  $k=1,...,3$. On the other hand, if $\bar{x}_1=\bar{x}_2=\bar{x}_3$ does not hold, then any dominated point is removed from the set \{$\bar{x}_1$, $\bar{x}_2$, $\bar{x}_3$\} (see Step 4(b(ii))). The latter case is typically encountered when the Pareto front and/or the domain is disconnected. Therefore, this algorithm is efficient in finding Pareto points even when the feasible set is discrete or disconnected. \\

\begin{description}
	\item[Step $\mathbf{1}$] { \textbf{(Input)}} \\ Choose the utopia point
	$u=(u_1,u_2,u_3)$ and $x \in X$. Set the
	number of partition points $N$ in the $(N+1)(N+2)/2$ grid points. Set $s=0$.
	\item[Step  $\mathbf{2}$] {\textbf{(Determine the individual minima)}}\\
	Solve Problem $(P_i):= \ds \min_{x\in  X} f_i(x), \;\; i=1,2,3$, that give the solutions $\bar{x}_{f_i}$, for $i=1,2,3$, respectively.\\
	Set $\bar{F}:= [f_1(\bar{x}_{f_3}),f_2(\bar{x}_{f_3}),f_3(\bar{x}_{f_3})]$
	
	\item[Step $\mathbf{3}$] { \textbf{(Generate grid points over the CHIM grid)}} \\
	Compute
	$a_k$, $b_k$ and $c_k$, for $k=1,2$, using \cite[Sections 5.1--5.2]{BurKayRiz2017}\,:  Let
	\[
	\varphi_k(\bar{x}_{f_i}) := \frac{\ds\prod_{j=1 \atop j\neq k}^3(f_j(\bar{x}_{f_i})-u_j)}
	{\ds\prod_{j=1,2}(f_j(\bar{x}_{f_i})-u_j) + \prod_{j=2,3}(f_j(\bar{x}_{f_i})-u_j) 
	+ \prod_{j=1,3}(f_j(\bar{x}_{f_i})-u_j)}\,,
	\]
	for $i = 1,2,3$.  Then
	\[
 	 a_k := \varphi_k(\bar{x}_{f_1})\,,\quad
 	 b_k := \varphi_k(\bar{x}_{f_2})\,,\quad 
 	 c_k := \varphi_k(\bar{x}_{f_3})\,.
	\]
	For $i_1:=0,1,...,(N-1)$\\
	\{\\
	\hspace*{8mm} Let $\hat{w}_k^{i_1}:=a_k+i_1(c_k-a_k)/N, \;\; k=1,2,$ \\
	\hspace*{8mm} and $\tilde{w}_r^{i_1}:=b_r+i_1(c_r-b_r)/N, \;\; r=1,2,$ \\
	\hspace*{8mm} For $i_2:=0,1,...,(N-i_1)$\\
	\hspace*{8mm} \{\\
	\hspace*{15mm} Let $w_\ell^{(i_1,i_2)}:=\hat{w}_\ell^{i_1}+i_2(\tilde{w}_\ell^{i_1}-\hat{w}_\ell^{i_1})/(N-i_1), \;\; \ell=1,2,$\\
	\hspace*{8mm} \}\\
	\}
	
	\item[Step $\mathbf{4}$] { \textbf{(Solve scalar sub-optimization problems)}}\\
	For $i_1:=0,1,...,(N-1)$\\
	\{\\
	\hspace*{8mm} For $i_2:=0,1,...,(N-i_1)$\\
	\hspace*{8mm} \{\\
	\hspace*{15mm} Let $w_k=w_k^{(i_1,i_2)},\;\; k=1,2.$ Let $w=(w_1,w_2,1-w_1-w_2)$.
	\begin{description}
		\item[(a)] Find $x_k$ that solves Problems ($P_w^k$), $k=1,2,3.$
		\item[(b)] Determine weak efficient points : \\
		Let $s:=s+1.$
		\begin{description}
			\item[(i)] If $\bar{x}_1=\bar{x}_2=\bar{x}_3$, then set $\bar{x}=\bar{x}_1$ (a weak efficient point)\\
			and  $F(s):=[f_1(\bar{x}),f_2(\bar{x}),f_3(\bar{x})].$
			\item [(ii)] If $\bar{x}_1=\bar{x}_2=\bar{x}_3$  does not hold, then, any dominated point is discarded by the following Steps~1-3.
			\begin{description}
				\item [(1)] If $f_2(\bar{x}_2) \geq f_2(\bar{x}_1)$ and $f_3(\bar{x}_3) \geq f_3(\bar{x}_1)$ then,  $\bar{x}_1$ weak efficient point,
				and $F(s):=[f_1(\bar{x}_1),f_2(\bar{x}_1),f_3(\bar{x}_1)].$
				\item [(2)] If $f_1(\bar{x}_1) \geq f_1(\bar{x}_2)$ and $f_3(\bar{x}_3) \geq f_3(\bar{x}_2)$ then,  $\bar{x}_2$ weak efficient point. Let $s:=s+1,$
				and $F(s):=[f_1(\bar{x}_2),f_2(\bar{x}_2),f_3(\bar{x}_2)].$
				\item [(3)] If $f_1(\bar{x}_1) \geq f_1(\bar{x}_3)$ and $f_2(\bar{x}_2) \geq f_2(\bar{x}_3)$ then,  $\bar{x}_3$ weak efficient point. Let $s:=s+1,$
				and $F(s):=[f_1(\bar{x}_3),f_2(\bar{x}_3),f_3(\bar{x}_3)].$
			\end{description}
			
		\end{description}
		\hspace*{8mm} \}\\
		\}
	\end{description}
	\item[Step $\mathbf{5}$] (Output)\\
	Set $F(s+1):= \bar{F}$.\\
	The array of Pareto points $F$, is an approximation of the Pareto front.
\end{description}

\subsection*{A.2\ \ Algorithm 7}

This algorithm is an extension of Algorithm~5  to the four-objective case.  In Algorithm~7, the Pareto front is approximated in six main steps. In Step 2, the individual minima of all four objectives are obtained. In Step~3, the boundaries of the Pareto front of the two-objective subproblems are constructed.  In Step~4, a linear programming problem is solved, to obtain base points that are used in the construction of the SBG grid in the spaces of the objective function triplets. In Step~5, the boundaries of the Pareto fronts of the three-objective subproblems are generated.  In Step~6, another linear programming problem is solved, to obtain base points that are used in the construction of the SBG grid in the space of all objective function quadruplets. Finally, in Step 7, a four-objective minimization problem is solved and so the points within the boundary of the Pareto front are generated.

\begin{description}
	\item[Step 1] (\textit{Input}) \\
	Set the number of partition points $N$.
	\item[Step 2] (\textit{Determine the individual minima}) \\
	Find $\bar{x}_{f_k}$ that solves Problem (P$_k$):=
	$\ds \min_{x\in X} f_k$, for $k = 1,2,3,4$. \\
	Set
	$F(k) := [f_1(\bar{x}_{f_k}), f_2(\bar{x}_{f_k}),
	f_3(\bar{x}_{f_k}),f_4(\bar{x}_{f_k})]$
	and $W(k):= e_k$, for $k = 1,2,3,4$.  \\
	Set $s = k+1$.
	\item[Step 3] (\textit{Generate Pareto fronts of the objective
		pairs}) \\
	\textit{Set the direction vector $v^{(i,k)}$ in the $f_if_k$-space,
		$i = 1,2,3$, $k = 2,3,4$, $i\neq k$}\,:\\
	Let $v^{(i,k)} := [v_1^{(i,k)},v_2^{(i,k)}]$. \\
	Set $v_1^{(i,k)} = \max\{f_i(\bar{x}_{f_i}),f_i(\bar{x}_{f_k})\} -
	f_i(\bar{x}_{f_i})$ \\
	and $v_2^{(i,k)} := \max\{f_k(\bar{x}_{f_i}),f_k(\bar{x}_{f_k})\} -
	f_k(\bar{x}_{f_k})$\,. \\
	{\em Auxiliary weights for a parallel ray in the $f_if_k$-space}\,: \\
	Let $c = 1/(v_1^{(i,k)} + v_2^{(i,k)})$,
	$\tilde{w}_1^{(i,k)} = cv_1^{(i,k)}$,
	$\tilde{w}_2^{(i,k)} = 1-\tilde{w}_1$. \\
	For $j:= 1,...,(N-1)$ \\
	\{ \hspace*{3mm}
	\begin{description}
		\item[(a)] (\textit{Generate Pareto points for $\{f_1,f_2\}$}) \\
		Let $w_1=(N-j)/N$; $w_2=1-w_1$;  $w_3=0$; $w_4=0$; \\ Set
		$W(s):=[w_1,w_2,w_3,w_4]$ and $s=s+1$.\\
		\textit{Form a base point for $\{f_1,f_2\}$}:\\
		Let $u_1 = w_1f_1(\bar{x}_{f_1}) + w_2f_1(\bar{x}_{f_2})$; $u_2 =
		w_1f_2(\bar{x}_{f_1}) + w_2f_2(\bar{x}_{f_2})$.
		\begin{description}
			\item[(i)]  Find $\bar{x}_1$ and $\bar{x}_2$ that solve
			Problems~(P$_w^1$) and (P$_w^2$), respectively.
			\item[(iii)] {\em Determine weak efficient points}\,:
			\begin{itemize}
				\item[$\bullet$] If $\bar{x}_1=\bar{x}_2=\bar{x}$ (a weak efficient point), then set \\
				$F(s):=[f_1(\bar{x}),f_2(\bar{x}),f_3(\bar{x}),f_4(\bar{x})]$ and $s=s+1$.
				\item[$\bullet$] If $\bar{x}_1\neq\bar{x}_2$, then $\bar{x}_1$ and $\bar{x}_2$ are weak efficient points if they are not dominated by each other. The following steps remove a dominated point, if it exists.
				\begin{description}
					\item [(1)] If $f_2(\bar{x}_2) \geq f_2(\bar{x}_1)$  then  $\bar{x}_1$ weak efficient point
					and \\
					 $F(s):=[f_1(\bar{x}_1),f_2(\bar{x}_1),f_3(\bar{x}_1),f_4(\bar{x}_1)].$
					\item [(2)] If $f_1(\bar{x}_1) \geq f_1(\bar{x}_2)$  then  $\bar{x}_2$ weak efficient point. Let $s:=s+1,$
					and $F(s):=[f_1(\bar{x}_2),f_2(\bar{x}_2),f_3(\bar{x}_2),f_4(\bar{x}_2)],$ and $W(s+1) = W(s)$.
				\end{description}
			\end{itemize}
		\end{description}
	\end{description}
	\begin{description}
		\item[(b)]  (\textit{Generate Pareto points for $\{f_2,f_3\}$}) \\
		Let $w_1 = 0$; $w_2 = (N-j)/N$; $w_3 = 1 - w_2$; $w_4 = 0$.  \\
		Do the same as in the rest of Step~3(a) for the new objective pair
		$\{f_2,f_3\}$ instead of $\{f_1,f_2\}$.
		\item[(c)]  (\textit{Generate Pareto points for $\{f_3,f_4\}$}) \\
		Let $w_1 = 0$; $w_2 = 0$; $w_3 = (N-j)/N$;  $w_4 = 1 - w_3$.  \\
		Do the same as in the rest of Step~3(a) for the new objective pair
		$\{f_3,f_4\}$ instead of $\{f_1,f_2\}$.
		\item[(d)]  (\textit{Generate Pareto points for $\{f_1,f_3\}$}) \\
		Let $w_1 = (N-j)/N$; $w_2 = 0$; $w_3 = 1 - w_1$; $w_4 = 0$.  \\
		Do the same as in the rest of Step~3(a) for the new objective pair
		$\{f_1,f_3\}$ instead of $\{f_1,f_2\}$.
		\item[(e)]  (\textit{Generate Pareto points for $\{f_1,f_4\}$}) \\
		Let  $w_1 = (N-j)/N$;  $w_2 = 0$; $w_3 = 0$; $w_4 = 1 - w_1$.  \\
		Do the same as in the rest of Step~3(a) for the new objective pair
		$\{f_1,f_4\}$ instead of $\{f_1,f_2\}$.
		\item[(f)]  (\textit{Generate Pareto points for $\{f_2,f_4\}$}) \\
		Let $w_1 = 0$; $w_2 = (N-j)/N$; $w_3 = 0$; $w_4 = 1 - w_2$.  \\
		Do the same as in the rest of Step~3(a) for the new objective pair
		$\{f_2,f_4\}$ instead of $\{f_1,f_2\}$.
	\end{description}
	\}
	\item[Step 4] (\textit{Form base points in the spaces of objective triplets}) \\
	Solve an LP problem as in \cite[Section	7.4]{MueGraeSch2009} to get an array of base points \\
	$U^{(1,2,3)}(r) = [u_1(r),u_2(r),u_3(r),u_4(r)]$, $r = 1,\ldots,(N-1)(N-2)/2$.\\
	Form base points $U^{(2,3,4)}(r)$, $U^{(1,3,4)}(r)$ and $U^{(1,2,4)}(r)$ in the $f_2f_3f_4$-, $f_1f_3f_4$- and $f_1f_2f_4$-spaces, respectively.
	
	\item[Step 5] (\textit{Generate Pareto fronts of the objective triplets})
	\begin{description}
		\item [(a)] (\textit{Generate interior Pareto points for $\{f_1,f_2,f_3\}$})\\
		\textit{Set the direction vector $v$ and weights for a parallel ray
			in the $f_1f_2f_3$-space}: \\
		Let $v := [v_1, v_2, v_3]$. Set $v_j := \max\{f_j(\bar{x}_{f_1}),
		f_j(\bar{x}_{f_2}), f_j(\bar{x}_{f_3})\} - f_j(\bar{x}_{f_j})$, $j=1,2,3$. \\
		Let $c=1/(v_1 +v_2 + v_3)$, $w_1= cv_1$, $w_2= cv_2$ and
		$w_3= 1-w_1-w_2$. \\
		For $r = 1,\ldots, \mbox{length}(U^{(1,2,3)}(r))$\\
		\{
		\begin{description}
			\item[(a)] Find $x_k$ that solves Problems ($P_w^k$), $k=1,2,3.$
			\item[(b)] Determine weak efficient points : \\
			Let $s:=s+1.$
			\begin{description}
				\item[(i)] If $\bar{x}_1=\bar{x}_2=\bar{x}_3$ then set, $\bar{x}=\bar{x}_1$ (a weak efficient point)\\
				and  $F(s):=[f_1(\bar{x}),f_2(\bar{x}),f_3(\bar{x}), f_4(\bar{x})].$
				\item [(ii)] If $\bar{x}_1=\bar{x}_2=\bar{x}_3$  does not hold, then, do the same as in Steps~4(b[ii]) of Algorithm~1.			
			\end{description}
		\end{description}
		\}
		\item [(b)] (\textit{Generate interior Pareto points for $\{f_2,f_3,f_4\}$, $\{f_1,f_3,f_4\}$ and $\{f_1,f_2,f_4\}$})\\
		Do the same as in the rest of Step~5(a) for the new objective triplets $\{f_2,f_3,f_4\}$, $\{f_1,f_3,f_4\}$ and $\{f_1,f_2,f_4\}$ instead of $\{f_1,f_2,f_3\}$, respectively.
		
	\end{description}
	\item[Step 6] (\textit{Form base points in the $f_1f_2f_3f_4$-space}) \\
	Given the Pareto points obtained in Steps~3 and 5 solve an LP problem
	to construct base points in the $f_1f_2f_3f_4$-space.
	\item[Step 7] (\textit{Generate interior Pareto points for
		$\{f_1,f_2,f_3,f_4\}$}) \\
	\textit{Set the direction vector $v$ and weights for a parallel ray
		in the $f_1f_2f_3f_4$-space}: \\
	Let $v := [v_1, v_2, v_3,v_4]$. Set $v_j := \max\{f_j(\bar{x}_{f_1}),
	f_j(\bar{x}_{f_2}), f_j(\bar{x}_{f_3}), f_j(\bar{x}_{f_4})\} - f_j(\bar{x}_{f_j})$, \\ $j=1,2,3,4$. \\
	Let $c=1/(v_1 +v_2 + v_3+v_4)$, $w_1= cv_1$, $w_2= cv_2$ , $w_3= cv_3$ and $w_4= 1-w_1-w_2-w_3$. \\
	For $r = 1,\ldots, \mbox{length}(U^{(1,2,3,4)}(r))$\\
	\{
	\begin{description}
		\item[(a)] Find $x_k$ that solves Problems ($P_w^k$), $k=1,2,3,4.$
		\item[(b)] Determine weak efficient points : \\
		Let $s:=s+1.$
		\begin{description}
			\item[(i)] If $\bar{x}_1=\bar{x}_2=\bar{x}_3=\bar{x}_4$ then set, $\bar{x}=\bar{x}_1$ (a weak efficient point)\\
			and  $F(s):=[f_1(\bar{x}),f_2(\bar{x}),f_3(\bar{x}), f_4(\bar{x})].$
			\item [(ii)] If $\bar{x}_1=\bar{x}_2=\bar{x}_3=\bar{x}_4$  does not hold, then, any dominated point is discarded by the following Steps~1-4.
			\begin{description}
				\item [(1)] If $f_2(\bar{x}_2) \geq f_2(\bar{x}_1)$ and $f_3(\bar{x}_3) \geq f_3(\bar{x}_1)$ and $f_4(\bar{x}_4) \geq f_4(\bar{x}_1)$ then,  $\bar{x}_1$ weak efficient point,
				and $F(s):=[f_1(\bar{x}_1),f_2(\bar{x}_1),f_3(\bar{x}_1), f_4(\bar{x}_1)].$
				\item [(2)] If $f_1(\bar{x}_1) \geq f_1(\bar{x}_2)$ and $f_3(\bar{x}_3) \geq f_3(\bar{x}_2)$  and $f_4(\bar{x}_4) \geq f_4(\bar{x}_2)$ then,  $\bar{x}_2$ weak efficient point. Let $s:=s+1,$
				and $F(s):=[f_1(\bar{x}_2),f_2(\bar{x}_2),f_3(\bar{x}_2),f_4(\bar{x}_2)].$
				\item [(3)] If $f_1(\bar{x}_1) \geq f_1(\bar{x}_3)$ and $f_2(\bar{x}_2) \geq f_2(\bar{x}_3)$ and $f_4(\bar{x}_4) \geq f_4(\bar{x}_3)$ then,  $\bar{x}_3$ weak efficient point. Let $s:=s+1,$
				and $F(s):=[f_1(\bar{x}_3),f_2(\bar{x}_3),f_3(\bar{x}_3),f_4(\bar{x}_3)].$
				\item [(4)] If $f_1(\bar{x}_1) \geq f_1(\bar{x}_4)$ and $f_2(\bar{x}_2) \geq f_2(\bar{x}_4)$ and $f_3(\bar{x}_3) \geq f_3(\bar{x}_4)$ then,  $\bar{x}_4$ weak efficient point. Let $s:=s+1,$
				and $F(s):=[f_1(\bar{x}_4),f_2(\bar{x}_4),f_3(\bar{x}_4),f_4(\bar{x}_4)].$
			\end{description}			
		\end{description}
	\end{description}
	\}
	\item[Step 8] (\textit{Output}) \\
	The array of Pareto points, $F$, is an approximation of the Pareto front.
\end{description}

\end{document}